\documentclass{amsart}
\usepackage{amssymb,latexsym,amsmath,amscd,graphicx,setspace,amsthm,verbatim,comment, stmaryrd,upgreek}
\usepackage{tikz,tikz-cd}
\usepackage[margin = 3 cm]{geometry} \input xy \xyoption{all}
\usepackage{tabularx}

\usetikzlibrary{shapes.geometric}

\usetikzlibrary{decorations.markings}

\tikzset{every loop/.style={min distance=10mm,looseness=10}}

\newcount\mycount

\theoremstyle{plain}
\newtheorem{theorem}{Theorem}[section]

\newtheorem{lthm}{Theorem} 

\makeatletter
\def\th@remark{%
  \thm@headfont{\bfseries}%
  \normalfont 
  \thm@preskip\topsep \divide\thm@preskip\tw@
  \thm@postskip\thm@preskip
}
\makeatother

\theoremstyle{remark} 
\newtheorem{remark}[theorem]{Remark}
      
\theoremstyle{definition}
\newtheorem{definition}[theorem]{Definition}

\begin{document} 
\title[On $\mathbb{Z}_{\ell}^{d}$-towers of graphs]{On $\mathbb{Z}_{\ell}^{d}$-towers of graphs}

\author{Sage DuBose, Daniel Valli\`{e}res}

\address{Mathematics and Statistics Department, California State University, Chico, CA 95929 USA}
\email{dvallieres@csuchico.edu}

\address{Mathematics and Statistics Department, California State University, Chico, CA 95929 USA}
\email{nsdubose@mail.csuchico.edu}

\date{\today}

\begin{abstract}
Let $\ell$ be a rational prime.  We show that an analogue of a conjecture of Greenberg in graph theory holds true.  More precisely, we show that when $n$ is sufficiently large, the $\ell$-adic valuation of the number of spanning trees at the $n$th layer of a $\mathbb{Z}_{\ell}^{d}$-tower of graphs is given by a polynomial in $\ell^{n}$ and $n$ with rational coefficients of total degree at most $d$ and of degree in $n$ at most one.
\end{abstract} 
\maketitle
\tableofcontents 
\tikzset{->-/.style={decoration={
  markings,
  mark=at position #1 with {\arrow{>}}},postaction={decorate}}}
\section{Introduction}
Throughout this paper, we fix a rational prime number $\ell$.  Let 
\begin{equation} \label{tower_one}
K = K_{0} \subseteq K_{1} \subseteq K_{2} \subseteq \ldots \subseteq K_{n} \subseteq \ldots 
\end{equation}
be a tower of number fields for which $K_{n}/K$ is Galois with ${\rm Gal}(K_{n}/K) \simeq \mathbb{Z}/\ell^{n}\mathbb{Z}$ for all $n \ge 0$.  A theorem of Iwasawa \cite[\S 4.2]{Iwasawa:1973} postulates the existence of non-negative integers $\mu,\lambda,n_{0}$ and an integer $\nu$ such that
$${\rm ord}_{\ell}(h_{n}) = \mu \ell^{n} + \lambda n + \nu,$$
for all $n \ge n_{0}$, where $h_{n}$ is the class number of $K_{n}$ and ${\rm ord}_{\ell}$ denotes the usual $\ell$-adic valuation on $\mathbb{Q}$.  The situation of constant field extensions for function fields over finite fields is classical and the situation of geometric $\mathbb{Z}_{\ell}$-extensions of global fields of characteristic $\ell$ was studied by Gold and Kisilevsky in \cite{Gold/Kisilevsky:1988}.  More recently, an analogous result was proven in the context of graph theory where field extensions are replaced with graph covers and the class number is replaced with the number of spanning trees.  See \cite{Vallieres:2021}, \cite{McGown/Vallieres:2022}, \cite{McGown/Vallieres:2022a} and \cite{Gonet:2021}, \cite{Gonet:2022}. 

Let $d$ be a fixed positive integer and consider a tower of number fields such as (\ref{tower_one}), where now $K_{n}/K$ is Galois with ${\rm Gal}(K_{n}/K) \simeq (\mathbb{Z}/\ell^{n}\mathbb{Z})^{d}$ for all $n \ge 0$.  Greenberg conjectured (see \cite[\S 7]{Cuoco/Monsky:1981}) that there exists a polynomial $P(X,Y) \in \mathbb{Q}[X,Y]$ of total degree at most $d$ and of degree at most $1$ in $Y$ such that
$${\rm ord}_{\ell}(h_{n}) = P(\ell^{n},n),$$
for $n$ large enough.  When $d=1$, the polynomial $P(X,Y)$ is simply Iwasawa's polynomial 
$$P(X,Y) = \mu X + \lambda Y + \nu \in \mathbb{Z}[X,Y].$$  
When $d\ge 2$, this conjecture of Greenberg was investigated in the number field situation by Cuoco and Monsky in \cite{Cuoco/Monsky:1981} and by Monsky in a series of papers (see \cite{Monsky:1981}, \cite{Monsky:1986}, \cite{Monsky:1986a}, \cite{Monsky:1987}, and \cite{Monsky:1989}).  More recently, Kleine gave in \cite{Kleine:2021} sufficient conditions in the number field situation when $d=2$ for Greenberg's conjecture to hold true.  Furthermore, for $d \ge 2$ the analogous conjecture in the function field case was proved by Wan in \cite{Wan:2019} for $\mathbb{Z}_{\ell}^{d}$-extensions of global fields of characteristic $\ell$ ramified at finitely many primes using the class number formula, $T$-adic $L$-functions and a result of Monsky (\cite[Theorem 5.6]{Monsky:1981}).  

The goal of the current paper is to adapt ideas of \cite{McGown/Vallieres:2022a} and \cite{Wan:2019} to prove the analogous conjecture in the context of graph theory.  In the category of finite graphs, one is lead to study sequences of covers of connected graphs 
\begin{equation} \label{gen_tower}
X = X_{0} \longleftarrow X_{1} \longleftarrow X_{2} \longleftarrow \ldots \longleftarrow X_{n} \longleftarrow \ldots
\end{equation}
for which the cover $X_{n}/X$ obtained by composing the covers $X_{n} \longrightarrow X_{n-1} \longrightarrow \ldots \longrightarrow X_{0}=X$ is Galois with group of covering transformations isomorphic to $(\mathbb{Z}/\ell^{n}\mathbb{Z})^{d}$.  (Throughout this paper, by a graph we shall mean what is sometimes referred to as a multigraph so that loops and parallel edges are allowed.)  Such a tower will be called a $\mathbb{Z}_{\ell}^{d}$-tower of graphs.  Our main result is the following:
\begin{lthm}[{Theorem~\ref{goal}}]\label{thmA}
Let $X$ be a connected graph and 
$$X = X_{0} \longleftarrow X_{1} \longleftarrow X_{2} \longleftarrow \ldots \longleftarrow X_{n} \longleftarrow \ldots$$
a $\mathbb{Z}_{\ell}^{d}$-tower of graphs.  Then, there exists $P(X,Y) \in \mathbb{Q}[X,Y]$ of total degree at most $d$ and of degree at most $1$ in $Y$ such that
$${\rm ord}_{\ell}(\kappa_{n}) = P(\ell^{n},n), $$
for $n$ large enough, where $\kappa_{n}$ is the number of spanning trees of $X_{n}$.
\end{lthm}

The paper is organized as follows.  In \S \ref{prel}, we gather together a few preliminaries on number theory, Pontryagin duals, graph theory and Artin-Ihara $L$-functions.  In \S \ref{voltage}, we remind the reader about voltage assignments in a way that is convenient for us.  We introduce in \S \ref{towers} the notion of a $\mathbb{Z}_{\ell}^{d}$-tower of graphs.  Following \cite{Wan:2019}, we remind the reader in \S \ref{char} about a useful character
$$\rho: \mathbb{Z}_{\ell}^{d} \longrightarrow \mathbb{Z}_{\ell}\llbracket T_{1},\ldots,T_{d}\rrbracket^{\times}$$ 
which we use in \S \ref{special} to study the special value at $u=1$ of Artin-Ihara $L$-functions as one goes up a $\mathbb{Z}_{\ell}^{d}$-tower of graphs.  This leads us to prove our main result, namely Theorem \ref{thmA}, in \S \ref{special}.  We end the paper with a few numerical examples in \S \ref{examples}.

\subsection*{Acknowledgement}
We would like to thank Randy Miller for lending us a computer to perform the numerical calculations presented in \S \ref{examples}.  DV would also like to thank C\'{e}dric Dion, Antonio Lei, Kevin McGown, and Anwesh Ray for several stimulating discussions on various topics related to the subject of this paper.

\section{Preliminaries} \label{prel}

\subsection{Number theory} 
Throughout this paper, we let $\mathbb{N} = \{0,1,2,\ldots\}$.  We need a few facts from the theory of cyclotomic number fields which can be found in \cite{Neukirch:1999} for instance.  The symbol $\overline{\mathbb{Q}}$ will denote a fixed algebraic closure of $\mathbb{Q}$ which we view as being embedded in $\mathbb{C}$ once and for all.  For a positive integer $m$, we let $\zeta_{m} = \exp(2\pi i/m) \in \overline{\mathbb{Q}}$ and we denote by $\mu_{m}$ the group of $m$th roots of unity.  We also fix a rational prime $\ell$ and the collection of all $\ell$-power roots of unity will be denoted by $\mu_{\ell^{\infty}}$.  For a positive integer $i$, there is a unique prime ideal $\mathcal{L}_{i}$ lying above $\ell$ in the ring of integers of the cyclotomic number field $\mathbb{Q}(\zeta_{\ell^{i}})$ and it is totally ramified.  From now on, we fix an embedding
\begin{equation} \label{embed}
\tau:\mathbb{Q}(\mu_{\ell^{\infty}}) \hookrightarrow \overline{\mathbb{Q}}_{\ell}, 
\end{equation}
and we view elements of $\mathbb{Q}(\mu_{\ell^{\infty}})$ inside of $\overline{\mathbb{Q}}_{\ell}$ via this embedding.  The extension of the $\ell$-adic valuation ${\rm ord}_{\ell}$ on $\mathbb{Q}$ to $\mathbb{C}_{\ell}$ will be denoted by $v_{\ell}$.   For $x \in \mathbb{Q}(\zeta_{\ell^{i}})$, one has
\begin{equation} \label{link_val}
v_{\ell}(\tau(x)) = \frac{1}{\varphi(\ell^{i})}{\rm ord}_{\mathcal{L}_{i}}(x),
\end{equation}
where $\varphi$ is the Euler $\varphi$-function, and ${\rm ord}_{\mathcal{L}_{i}}$ is the valuation on $\mathbb{Q}(\zeta_{\ell^{i}})$ associated to the prime ideal $\mathcal{L}_{i}$.  The absolute value on $\mathbb{C}_{\ell}$ extending the $\ell$-adic absolute value on $\mathbb{Q}$ normalized in the usual way will be denoted by $|\cdot|_{\ell}$.  Note that it follows from (\ref{link_val}) that 
\begin{equation} \label{in_disk}
|1-\xi|_{\ell} < 1
\end{equation}
for all $\xi \in \mu_{\ell^{\infty}}$, since $\mathcal{L}_{i} = (1 - \zeta_{\ell^{i}})$.

\subsection{Pontryagin dual}
For the basic facts contained in this section, we refer the reader to the first chapter of \cite{Neukirch/Schmidt/Wingberg:2008}.  If $G$ is an abelian compact topological group, we let $G^{\vee}$ denote the Pontryagin dual of $G$.  It consists of the group of continuous group morphisms from $G$ to $S^{1}$, where $S^{1}$ denotes the unit circle in $\mathbb{C}^{\times}$ with its induced topology.  An element of $G^{\vee}$ will simply be called a character of $G$.  If $G$ is an abelian profinite group, then every character of $G$ has finite image and thus takes values in $\mu_{\infty}$, where $\mu_{\infty}$ is the collection of all roots of unity in $S^{1}$.  It follows that if $G$ is an abelian profinite group, then
\begin{equation} \label{disc_top}
G^{\vee} = G^{*}, 
\end{equation}
where $G^{*}$ is the group of continuous group morphisms from $G$ to $\mu_{\infty}$, where $\mu_{\infty}$ is induced with the discrete topology.  Now, if $G$ is an abelian pro-$\ell$ group, then any character actually takes values in $\mu_{\ell^{\infty}}$.  Thus, via the embedding (\ref{embed}), we will view the characters of an abelian pro-$\ell$ group as taking values in $\overline{\mathbb{Q}}_{\ell}$.  Furthermore, it follows from (\ref{disc_top}) that
\begin{equation} \label{cont_char}
G^{\vee} \subseteq {\rm Hom}_{{\rm cont}}(G,\mathbb{C}_{\ell}^{\times}).
\end{equation}
The group $G$ arising in this paper will be $\mathbb{Z}_{\ell}^{d}$, where $d$ is a fixed positive integer.  For $i=1,\ldots,d$, let $G_{i} = \mathbb{Z}_{\ell}$, so that 
$$\mathbb{Z}_{\ell}^{d} = G_{1} \times \ldots \times G_{d} \text{ and } (\mathbb{Z}_{\ell}^{d})^{\vee} \simeq G_{1}^{\vee} \times \ldots \times G_{d}^{\vee}. $$
If $\psi_{i} \in G_{i}^{\vee}$ for $i=1,\ldots,d$, then the character of $(\mathbb{Z}_{\ell}^{d})^{\vee}$ corresponding to $(\psi_{1},\ldots,\psi_{d}) \in G_{1}^{\vee}\times \ldots \times G_{d}^{\vee}$ via the second isomorphism above will be denoted by $\psi_{1}\otimes \ldots \otimes \psi_{d}$.  Thus, if $a = (a_{1},\ldots,a_{d}) \in \mathbb{Z}_{\ell}^{d}$ and $\psi \in (\mathbb{Z}_{\ell}^{d})^{\vee}$, then 
$$\psi(a) = \psi_{1} \otimes \ldots \otimes \psi_{d}(a_{1},\ldots,a_{d}) = \psi_{1}(a_{1}) \cdot \ldots \cdot \psi_{d}(a_{d}), $$
where $\psi_{i}$ is the restriction of $\psi$ to $G_{i}$.

\subsection{Graph theory}
Our main references for this section are \cite{Serre:1977} and \cite{Sunada:2013}.  A graph $X = (V_{X},\mathbf{E}_{X})$ thus consists of a set of vertices $V_{X}$, a set of directed edges $\mathbf{E}_{X}$, an incidence map ${\rm inc}: \mathbf{E}_{X} \longrightarrow V_{X} \times V_{X}$ denoted by $e \mapsto {\rm inc}(e) = (o(e),t(e))$, and an inversion map $\mathbf{E}_{X} \longrightarrow \mathbf{E}_{X}$ denoted by $e \mapsto \bar{e}$ satisfying $\bar{\bar{e}}=e$, $\bar{e} \neq e$, and $o(e) = t(\bar{e})$ for all $e \in \mathbf{E}_{X}$.  The set of undirected edges will be denoted by $E_{X}$.  Throughout this paper, we allow loops and multiple edges.  Such graphs are sometimes called multigraphs, but we will refer to them more simply as graphs.  A graph that has neither loops nor multiple edges will be called a simple graph.  \emph{All graphs in this paper will be assumed to be finite meaning that both $V_{X}$ and $\mathbf{E}_{X}$ are finite sets.}  If $v \in V_{X}$, then we let
$$\mathbf{E}_{X,v} = \{e \in \mathbf{E}_{X} \, | \, o(e) = v \} \text{ and } {\rm val}_{X}(v) = |\mathbf{E}_{X,v}|. $$

Let $Y= (V_{Y},\mathbf{E}_{Y})$ and $X = (V_{X},\mathbf{E}_{X})$ be two graphs.  Recall from \cite{Sunada:2013} that a morphism of graphs $f = (f_{V},f_{E}):Y \longrightarrow X$ consists of two functions $f_{V}:V_{Y} \longrightarrow V_{X}$ and $f_{E}:\mathbf{E}_{Y} \longrightarrow \mathbf{E}_{X}$ satisfying
\begin{enumerate}
\item $f_{V}(o(e)) = o(f_{E}(e))$ for all $e \in \mathbf{E}_{Y}$,
\item $f_{V}(t(e)) = t(f_{E}(e))$ for all $e \in \mathbf{E}_{Y}$,
\item $\overline{f_{E}(e)} = f_{E}(\bar{e})$ for all $e \in \mathbf{E}_{Y}$.
\end{enumerate}
A morphism is an isomorphism if and only if both $f_{V}$ and $f_{E}$ are bijections.  The group of automorphisms of a graph $X$ will be denoted by ${\rm Aut}(X)$.  We will often drop the indices $V$ and $E$ for the functions $f_{V}$ and $f_{E}$ and denote them both simply by $f$.

From now on, let $Y$ and $X$ be two connected graphs.  A morphism of graphs $\pi:Y \longrightarrow X$ is called a covering map (or just a cover) if the following two extra conditions are satisfied:
\begin{enumerate}
\item $\pi:V_{Y} \longrightarrow V_{X}$ is surjective,
\item For all $w \in V_{Y}$, the restriction of $\pi$ to $\mathbf{E}_{Y,w}$, denoted by $\pi|_{\mathbf{E}_{Y,w}}$ induces a bijection
$$\pi|_{\mathbf{E}_{Y,w}}:  \mathbf{E}_{Y,w} \stackrel{\approx}{\longrightarrow} \mathbf{E}_{X,\pi(w)}. $$
\end{enumerate}
If $\pi:Y\longrightarrow X$ is a covering map, then we let
$${\rm Aut}_{\pi}(Y/X) = \{ \phi \in {\rm Aut}(Y) \, : \, \pi \circ \phi = \pi \}, $$
and we shall often drop $\pi$ from the notation if it is clear from the context.  Clearly, ${\rm Aut}(Y/X)$ is a subgroup of ${\rm Aut}(Y)$.  A cover $\pi:Y \longrightarrow X$ is called Galois if the action of ${\rm Aut}(Y/X)$ on the fiber $\pi^{-1}(v)$ for all $v \in V_{X}$ is transitive.  In this case, we write ${\rm Gal}(Y/X)$ instead of ${\rm Aut}(Y/X)$.  If $Y/X$ is a Galois cover, then we have the usual Galois correspondence between subgroups of ${\rm Gal}(Y/X)$ and equivalence classes of intermediate covers of $Y/X$.  By an abelian cover of graphs $Y/X$, we mean a Galois covering map $\pi:Y \longrightarrow X$ for which ${\rm Gal}(Y/X)$ is an abelian group.

\subsection{Ihara zeta functions and Artin-Ihara $L$-functions}
For more details regarding this section, we refer the reader to \cite{Terras:2011}.  Let $X$ be a connected graph and label the vertices $V_{X} = \{v_{1},\ldots,v_{g} \}$.  \emph{From now on, we assume that all of our graphs do not have vertices of valency one.}  We denote the Ihara zeta function of $X$ by $Z_{X}(u)$.  (This zeta function was first defined by Ihara in \cite{Ihara:1966} and interpreted in terms of graph theory by Sunada in \cite{Sunada:1986} following a suggestion of Serre.)  The three-term determinant formula \cite[Theorem 2.5]{Terras:2011} gives
$$Z_{X}(u)^{-1} = (1-u^{2})^{-\chi(X)} {\rm det}(I - Au + (D-I)u^{2}), $$
where $I$ is the identity matrix, $A$ the adjacency matrix of $X$, $D$ the degree (or valency) matrix of $X$, and $\chi(X)$ the Euler characteristic of $X$.  From now on, we let
$$h_{X}(u) = {\rm det}(I - Au + (D-I)u^{2}) \in \mathbb{Z}[u].$$
Furthermore, one has $h_{X}(1) = 0$, since the Laplacian matrix $D-A$ is singular, and the main result of \cite{Northshield:1998} and also \cite[Theorem B]{Hashimoto:1990} imply
\begin{equation} \label{class_formula}
h_{X}'(1) = -2 \chi(X) \kappa_{X},
\end{equation}
where $\kappa_{X}$ is the number of spanning trees of $X$.  

More generally, let $Y/X$ be an abelian cover of connected graphs and let $G = {\rm Gal}(Y/X)$.  (The theory can be extended to arbitrary Galois covers not necessarily abelian.  See \cite{Terras:2011} and also \cite{Stark/Terras:1996}, \cite{Stark/Terras:2000}.)  If $\psi \in G^{\vee}$, then we denote the corresponding Artin-Ihara $L$-function by $L_{X}(u,\psi)$.  For each $i=1,\ldots,g$, we let $w_{i} \in V_{Y}$ be a fixed vertex in the fiber of $v_{i}$.  The three-term determinant formula \cite[Theorem 18.15]{Terras:2011} gives
$$L_{X}(u,\psi)^{-1} = (1-u^{2})^{-\chi(X)} {\rm det}(I - A_{\psi}u + (D-I)u^{2}), $$
where $I$ is the identity matrix, $A_{\psi}$ the twisted adjacency matrix associated to $\psi$, $D$ the degree matrix of $X$, and $\chi(X)$ the Euler characteristic of $X$.  Let us remind the reader how $A_{\psi}$ is defined (see \cite[Definition 18.13]{Terras:2011}).  For $\sigma \in G$, one sets $A(\sigma)=(a_{ij}(\sigma))$ to be the $g \times g$ matrix defined via
\begin{equation*}
a_{ij}(\sigma) = 
\begin{cases}
\text{Twice the number of undirected loops at the vertex }w_{i}, &\text{ if } i=j \text{ and } \sigma = 1;\\
\text{The number of undirected edges connecting $w_{i}$ to $w_{j}^{\sigma}$}, &\text{ otherwise},
\end{cases}
\end{equation*}
and if $\psi \in G^{\vee}$, then the twisted adjacency matrix $A_{\psi}$ is defined via
$$A_{\psi} = \sum_{\sigma \in G} \psi(\sigma) \cdot A(\sigma). $$
From now on, we let
\begin{equation} \label{3_term}
h_{X}(u,\psi) =  {\rm det}(I - A_{\psi}u + (D-I)u^{2}) \in \mathbb{Z}[\psi][u]. 
\end{equation}
Throughout this paper, we let $\psi_{0}$ denote the trivial character.  We have $A_{\psi_{0}} = A$, the usual adjacency matrix of $X$, so that 
$$Z_{X}(u) = L_{X}(u,\psi_{0}) \text{ and } h_{X}(u) = h_{X}(u,\psi_{0}). $$
The Artin-Ihara $L$-functions satisfy the usual Artin formalism (see \cite[Proposition 18.10]{Terras:2011}), and thus using the fact that
\begin{equation} \label{Riemann-Hurwitz}
\chi(Y) = |G| \cdot \chi(X),
\end{equation}
one has
\begin{equation} \label{prod}
h_{Y}(u) = h_{X}(u) \prod_{\psi \neq \psi_{0}} h_{X}(u,\psi),
\end{equation}
where the product is over all non-trivial characters $\psi$ of $G$.  Differentiating (\ref{prod}), evaluating at $u=1$, and using (\ref{class_formula}) and (\ref{Riemann-Hurwitz}) leads to the formula
\begin{equation} \label{imp_formula}
|G|\cdot \kappa_{Y} = \kappa_{X} \prod_{\psi \neq \psi_{0}}h_{X}(1,\psi),
\end{equation}
assuming that $\chi(X) \neq 0$.  \emph{From now on, we assume that $\chi(X) \neq 0$ throughout.}

\section{Construction of Galois covers via voltage assignments} \label{voltage}
A convenient way to construct Galois covers of graphs is via voltage assignments.  Our main reference for this section is \cite{Gross/Tucker:2001}.  (See also \cite{Gross/Tucker:1977}.)  Let $X$ be a graph and let $\gamma:E_{X} \longrightarrow \mathbf{E}_{X}$ be a section of the natural map $\mathbf{E}_{X} \longrightarrow E_{X}$.  We set $S = \gamma(E_{X})$.  If $G$ is a finite abelian group (the abelian condition could be removed if desired) and $\alpha:S \longrightarrow G$ is a function (sometimes called a voltage assignment), then we extend $\alpha$ to $\bar{S}$ by setting
$$\alpha(\bar{s}) = \alpha(s)^{-1}, $$
when $s \in S$.  The derived graph $X(G,S,\alpha)$ is a new graph constructed as follows.  The set of vertices $V$ is $V_{X} \times G$ and the set of directed edges $\mathbf{E}$ is $\mathbf{E}_{X} \times G$.  The incidence map is given by
$${\rm inc}((e,\sigma)) = ((o(e),\sigma),(t(e),\sigma \cdot \alpha(e))),$$
and the inversion map by
$$\overline{(e,\sigma)} = (\bar{e},\sigma \cdot \alpha(e)). $$
We leave it to the reader to check that $X(G,S,\alpha)$ is a graph.  

Furthermore, if for simplicity we let $Y = X(G,S,\alpha)$ and we define $p:Y \longrightarrow X$ via
\begin{equation} \label{cover_map_volt}
p((v,\sigma)) = v \text{ and } p((e,\sigma)) = e, 
\end{equation}
then $p$ is a morphism of graphs.  Under the assumption that both $X$ and $Y$ are connected, the map $p$ is in fact a Galois cover.  For each $\sigma \in G$, we let $\varphi_{\sigma} \in {\rm Aut}(Y)$ be defined via
$$\varphi_{\sigma}((v,\tau)) = (v, \sigma \cdot \tau) \text{ and } \varphi_{\sigma}((e,\tau))= (e,\sigma \cdot \tau).$$
Then the function $\phi:G \longrightarrow {\rm Aut}_{p}(Y/X)$ given by $\sigma \mapsto \phi(\sigma) = \varphi_{\sigma}$ can be checked to be a group isomorphism.  In other words, the cover $p:Y\longrightarrow X$ is a Galois cover whose group of covering transformations is isomorphic to $G$.  

Any abelian cover arises from a derived graph as above as we now explain.  If $\pi:Y\longrightarrow X$ is an abelian cover of connected graphs with Galois group $G$, then consider the natural surjective group morphism
$$\mu:H_{1}(X,\mathbb{Z}) \longrightarrow {\rm Gal}(Y/X), $$
and let $T$ be a spanning tree of $X$.  If $e \in \mathbf{E}_{X}$, then let $c_{e} = d_{e} \cdot e$, where $d_{e}$ is the unique geodesic path going from $t(e)$ to $o(e)$ within $T$.  Then $c_{e}$ is a closed path in $X$.  Let $S$ be the image of a section $\gamma$ as above, and define $\alpha:S \longrightarrow G$ via $s \mapsto \alpha(s) = \mu(\langle c_{s}\rangle)$, where $\langle c_{s} \rangle \in H_{1}(X,\mathbb{Z})$ is the corresponding cycle.  Consider now the graph $X(G,S,\alpha)$ and choose a labeling $V_{X} = \{v_{1},\ldots,v_{g} \}$.  Let $\widetilde{T}$ be a lift of $T$ to $Y$, and set $V_{\widetilde{T}} = \{w_{1},\ldots,w_{g}\}$, where $\pi(w_{i}) = v_{i}$ for $i=1,\ldots,g$.  Define $\phi:X(G,S,\alpha) \longrightarrow Y$ via
\begin{equation} \label{iso_con}
\phi((v_{i},\sigma)) = \sigma(w_{i}) \text{ and }  \phi((e,\sigma)) = \sigma(\tilde{e}),
\end{equation}
where $\tilde{e}$ is the unique lift of $e$ to $Y$ starting at $w_{i}$ if $o(e)=v_{i}$.  Then, one can check that $\phi$ is an isomorphism of graphs and that the following diagram 
\begin{equation*}
\begin{tikzcd}
X(G,S,\alpha) \arrow["\phi",rr] \arrow["p",rd]  & & Y \arrow["\pi",ld] \\
& X&\\ 
\end{tikzcd}
\end{equation*}
commutes.

Now, let $G_{1}$ and $G_{2}$ be both finite abelian groups, and let $f:G_{1} \longrightarrow G_{2}$ be a group morphism.  If we start with a function $\alpha:S \longrightarrow G_{1}$, then we can also consider the function $f\circ \alpha:S \longrightarrow G_{2}$ and the graphs $Y_{1} = X(G_{1},S,\alpha)$ and $Y_{2} = X(G_{2},S,f \circ \alpha)$.  We leave it to the reader to check that $f_{*}:Y_{1} \longrightarrow Y_{2}$ defined via
$$f_{*}(v,\sigma_{1}) = (v,f(\sigma_{1})) \text{ and } f_{*}(e,\sigma_{1}) = (e,f(\sigma_{1})) $$
is a morphism of graphs.  Moreover, if $f$ is surjective and both $Y_{1}$ and $Y_{2}$ are connected, then $f_{*}$ is a Galois cover with Galois group isomorphic to ${\rm ker}(f)$.  Letting $p_{i}$ be the covering map $p_{i}:Y_{i} \longrightarrow X$ defined in (\ref{cover_map_volt}) above, we obtain that $(Y_{2},f_{*},p_{2})$ is an intermediate cover of $p_{1}:Y_{1} \longrightarrow X$, in other words the following diagram
\begin{equation*}
\begin{tikzcd}
X(G_{1},S,\alpha) \arrow["p_{1}",dd] \arrow[rd,"f_{*}"] & \\
& X(G_{2},S,f \circ \alpha) \arrow[ld,"p_{2}"]\\
X & 
\end{tikzcd}
\end{equation*}
commutes.

\section{$\mathbb{Z}_{\ell}^{d}$-towers of graphs} \label{towers}
Let $X$ be a connected graph and $d$ a fixed positive integer.  We start with the following definition.
\begin{definition} \label{d_tower}
A $\mathbb{Z}_{\ell}^{d}$-tower of graphs above a connected graph $X$ consists of a sequence of covers of connected graphs
$$X = X_{0} \longleftarrow X_{1} \longleftarrow X_{2} \longleftarrow \ldots \longleftarrow X_{n} \longleftarrow \ldots $$
having the property that the cover $X_{n}/X$ obtained from the composition $X_{n} \longrightarrow \ldots \longrightarrow X_{2} \longrightarrow X_{1} \longrightarrow X$ is Galois with group of covering transformations isomorphic to $(\mathbb{Z}/\ell^{n}\mathbb{Z})^{d}$.  
\end{definition}
Note that when $d=1$, such towers were simply called abelian $\ell$-towers in \cite{Vallieres:2021}, \cite{McGown/Vallieres:2022}, \cite{McGown/Vallieres:2022a} and \cite{Lei/Vallieres:2022}.  Let now $S$ be the image of a section of the natural map $\mathbf{E}_{X} \longrightarrow E_{X}$ as explained in \S \ref{voltage}.  If we start with a function $\alpha:S \longrightarrow \mathbb{Z}_{\ell}^{d}$, then for each $n \in \mathbb{N}$ we let $\alpha_{n}:S \longrightarrow (\mathbb{Z}/\ell^{n}\mathbb{Z})^{d}$ be the function obtained from the compositions
$$\mathbb{Z}_{\ell}^{d} \longrightarrow (\mathbb{Z}_{\ell}/\ell^{n}\mathbb{Z}_{\ell})^{d} \stackrel{\simeq}{\longrightarrow} (\mathbb{Z}/\ell^{n}\mathbb{Z})^{d}. $$
To simplify the notation, we now let $G(n) = (\mathbb{Z}/\ell^{n}\mathbb{Z})^{d}$.  Under the assumption that all graphs $X(G(n),S,\alpha_{n})$ are connected, it follows from \S \ref{voltage} that we obtain a $\mathbb{Z}_{\ell}^{d}$-tower of graphs
\begin{equation} \label{concrete_tower}
X \longleftarrow X(G(1),S,\alpha_{1}) \longleftarrow X(G(2),S,\alpha_{2}) \longleftarrow \ldots \longleftarrow X(G(n),S,\alpha_{n}) \longleftarrow \ldots 
\end{equation}

Conversely every $\mathbb{Z}_{\ell}^{d}$-tower of graphs in the sense of Definition \ref{d_tower} is isomorphic in a suitable sense to a $\mathbb{Z}_{\ell}^{d}$-tower as in (\ref{concrete_tower}) for some function $\alpha:S \longrightarrow \mathbb{Z}_{\ell}^{d}$ as we now briefly explain.  Let 
$$X \longleftarrow X_{1} \longleftarrow X_{2} \longleftarrow \ldots \longleftarrow X_{n} \longleftarrow \ldots $$
be a $\mathbb{Z}_{\ell}^{d}$-tower of graphs and let $G_{n} = {\rm Gal}(X_{n}/X)$.  Then we have group morphisms
$$\mu_{n}:H_{1}(X,\mathbb{Z}) \longrightarrow G_{n} $$
that are compatible so that we get a group morphism
$$\mu:H_{1}(X,\mathbb{Z}) \longrightarrow \lim_{\substack{\longleftarrow \\ n \ge1 }}G_{n} \simeq \mathbb{Z}_{\ell}^{d}. $$
Let $T$ be a spanning tree of $X$, and define $\alpha:S \longrightarrow \mathbb{Z}_{\ell}^{d}$ via $s \mapsto \alpha(s) = \mu(\langle c_{s} \rangle)$, where $c_{s} = d_{s} \cdot s$ and $d_{s}$ is the unique geodesic path going from $t(s)$ to $o(s)$ within $T$.  Let now $T_{1}$ be a lift of $T$ to $X_{1}$ via the covering map $X_{1} \longrightarrow X$.  Then, we get the isomorphism of graphs
$$\phi_{1}:X(G(1),S,\alpha_{1}) \longrightarrow X_{1}$$
as defined in (\ref{iso_con}) above.  Let now $T_{2}$ be a lift of $T_{1}$ to $X_{2}$ via the covering map $X_{2} \longrightarrow X_{1}$.  Then, $T_{2}$ is a lift of $T$ to $X_{2}$, and we get again from (\ref{iso_con}) an isomorphism of graphs
$$\phi_{2}:X(G(2),S,\alpha_{2}) \longrightarrow X_{2}.$$
Keeping going like this, one constructs isomorphisms of graphs $\phi_{n}:X(G(n),S,\alpha_{n}) \longrightarrow X_{n}$ for which all the squares and the triangle in the diagram
\begin{equation} \label{nice_dia}
\begin{tikzcd}
         & \arrow[ld] X(G(1),S,\alpha_{1})  \arrow[dd, "\phi_{1}"]  & \arrow[l]X(G(2),S,\alpha_{2})   \arrow[dd, "\phi_{2}"] &\arrow[l] \ldots  &\arrow[l] \arrow[dd, "\phi_{n}"]  X(G(n),S,\alpha_{n}) & \arrow[l]  \ldots\\  
         X &  & & & & \\
          & \arrow[lu]X_{1}  & \arrow[l]X_{2}   & \arrow[l] \ldots  & \arrow[l] X_{n}  \arrow[l] & \arrow[l] \ldots
\end{tikzcd}
\end{equation}
commute.

\section{An $(\ell,T_{1},\ldots,T_{d})$-adic character} \label{char}
Throughout this paper, for a fixed positive integer $d$ and a unital commutative ring $R$, we set 
$$R[T] = R[T_{1},\ldots,T_{d}] \text{ and }  R\llbracket T \rrbracket = R\llbracket T_{1},\ldots,T_{d} \rrbracket.$$  
The ring $\mathbb{Z}_{\ell}\llbracket T \rrbracket$ is a local ring with maximal ideal $\mathfrak{m} = (\ell,T_{1},\ldots,T_{d})$, and we consider it with its usual $\mathfrak{m}$-adic topology.  For any $a \in \mathbb{Z}_{\ell}$ and any $i=1,\ldots,d$, the function $\mathbb{Z}_{\ell} \longrightarrow \mathbb{Z}_{\ell}\llbracket T_{i} \rrbracket \subseteq \mathbb{Z}_{\ell}\llbracket T \rrbracket$ defined by $a \mapsto (1-T_{i})^{a}$ is continuous.  (See for instance \cite[Proposition 3.1]{McGown/Vallieres:2022a} and the discussion that follows.)  Therefore, one gets a continuous group morphism $\rho:\mathbb{Z}_{\ell}^{d} \longrightarrow \mathbb{Z}_{\ell}\llbracket T \rrbracket^{\times}$ defined via
$$a = (a_{1},\ldots,a_{d}) \mapsto \rho(a) = (1-T_{1})^{a_{1}}\cdot \ldots \cdot (1-T_{d})^{a_{d}}. $$
We set 
$$Q_{a}(T) := \rho(a), $$
and note that if $a \in \mathbb{N}^{d}$, then $Q_{a}(T) \in \mathbb{Z}[T]$ and if $a \in \mathbb{Z}^{d}$, then $Q_{a}(T) \in \mathbb{Z}\llbracket T \rrbracket$.  

From now on, we let 
$$D = \{t = (t_{1},\ldots,t_{d}) \in \mathbb{C}_{\ell}^{d} \, : \, |t_{i}|_{\ell} < 1 \text{ for all } i=1,\ldots,d \}. $$
Given any $t \in D$, the function ${\rm ev}_{t}:\mathbb{Z}_{\ell}\llbracket T \rrbracket \longrightarrow \mathbb{C}_{\ell}$ given by
\begin{equation} \label{cont_ev}
Q(T) \mapsto {\rm ev}_{t}(Q(T)) = Q(t) 
\end{equation}
is continuous.  Indeed, if $Q(T) \in \mathfrak{m}^{N}$ for some positive integer $N$, and if we write
$$Q(T) = \sum_{i=0}^{\infty}Q_{i}(T), $$
where $Q_{i}(T) \in \mathbb{Z}_{\ell}[T]$ is homogeneous of degree $i$, then
\begin{equation*}
\begin{aligned}
|Q(t)|_{\ell} &\le \sum_{i=0}^{N-1}|Q_{i}(t)|_{\ell} + \sum_{i = N}^{\infty}|Q_{i}(t)|_{\ell}\\
&\le \frac{1}{\ell^{N}}\binom{N-1 + d -1}{d-1}\sum_{i=0}^{N-1}(\ell x)^{i} + \sum_{\substack{\alpha \in \mathbb{N}^{d} \\ |\alpha| \ge N}}|t^{\alpha}|_{\ell},
\end{aligned}
\end{equation*}
where we use the usual multi-index notation for the sum on the right and $x = {\rm max}\{ |t_{1}|_{\ell},\ldots,|t_{d}|_{\ell}\}$.  Now, both these sums can be made arbitrarily small for $N$ large, and this shows that (\ref{cont_ev}) is continuous.

Furthermore, if $\psi \in (\mathbb{Z}_{\ell}^{d})^{\vee}$, then the function $\mathbb{Z}_{\ell}^{d} \longrightarrow \mathbb{C}_{\ell}$ defined via 
$$a \mapsto \psi(a) $$
is continuous as well by (\ref{cont_char}).

Following \cite{Wan:2019}, for $\psi \in (\mathbb{Z}_{\ell}^{d})^{\vee}$, we let
$$t_{\psi} = (1 - \psi_{1}(1),\ldots,1-\psi_{d}(1)) \in D, $$
where $\psi = \psi_{1} \otimes \psi_{2} \otimes \ldots \otimes \psi_{d}$.  Note that $t_{\psi} \in D$ by (\ref{in_disk}).  It follows from the above discussion that for a fixed $\psi \in (\mathbb{Z}_{\ell}^{d})^{\vee}$ both functions ${\rm ev}_{t_{\psi}} \circ \rho$ and $\psi$ are continuous functions from $\mathbb{Z}_{\ell}^{d}$ to $\mathbb{C}_{\ell}$. Since they are equal on $\mathbb{N}^{d}$, which is dense in $\mathbb{Z}_{\ell}^{d}$, one gets
\begin{equation} \label{imp_equal}
{\rm ev}_{t_{\psi}} \circ \rho(a) = Q_{a}(t_{\psi}) = \psi(a),
\end{equation}
for all $a \in \mathbb{Z}_{\ell}^{d}$.

\begin{remark}
In \cite{Wan:2019}, the author is working with the character
$$a = (a_{1},\ldots,a_{d}) \mapsto (1+T_{1})^{a_{1}} \cdot \ldots \cdot (1+T_{d})^{a_{d}} $$
and the classical points
$$t_{\psi} = (\psi_{1}(1)-1,\ldots,\psi_{d}(1)-1). $$
We choose to work with our current convention in order to align with the previous work contained in \cite{Vallieres:2021}, \cite{McGown/Vallieres:2022}, and \cite{McGown/Vallieres:2022a}.
\end{remark}

\section{The special value at $u=1$ of Artin-Ihara $L$-functions in $\mathbb{Z}_{\ell}^{d}$-towers of graphs} \label{special}

Let $X$ be a graph (with loops and multiple edges allowed) and recall that we are assuming $X$ is finite, connected, with no vertex of degree one, and $\chi(X) \neq 0$.  Fix also $S$ as explained in \S \ref{voltage}.  We start with the following theorem.
\begin{theorem} \label{main_thm}
Let $X$ be as above, $\alpha:S \longrightarrow \mathbb{Z}_{\ell}^{d}$ a function for which all the graphs $X(G(n),S,\alpha_{n})$ are connected, and label the vertices $V_{X} = \{v_{1},\ldots,v_{g} \}$.  Consider the $\mathbb{Z}_{\ell}^{d}$-tower of graphs
$$X \longleftarrow X(G(1),S,\alpha_{1}) \longleftarrow X(G(2),S,\alpha_{2}) \longleftarrow \ldots \longleftarrow X(G(n),S,\alpha_{n}) \longleftarrow \ldots $$
and let
$$Q(T) = {\rm det}(D - A_{\rho}) \in \mathbb{Z}_{\ell}\llbracket T\rrbracket, $$
where $D$ is the degree matrix of $X$ and
\begin{equation*}
A_{\rho} = \left(\sum_{\substack{s \in S \\ {\rm inc}(s) = (v_{i},v_{j})}}\rho(\alpha(s)) + \sum_{\substack{s \in S \\ {\rm inc}(s) = (v_{j},v_{i})}}\rho(-\alpha(s)) \right) \in M_{g \times g}(\mathbb{Z}_{\ell}\llbracket T \rrbracket).
\end{equation*}
Then, for all $n \in \mathbb{N}$ and for all $\psi \in G(n)^{\vee}$, one has
$$Q(t_{\widetilde{\psi}}) = h_{X}(1,\psi),$$
where $\widetilde{\psi} \in (\mathbb{Z}_{\ell}^{d})^{\vee}$ is obtained from $\psi$ after composing with the natural projection map $\mathbb{Z}_{\ell}^{d} \twoheadrightarrow G(n)$.  
\end{theorem}
\begin{proof}
Given $t \in D$, the function ${\rm ev}_{t}:\mathbb{Z}_{\ell}\llbracket T\rrbracket \longrightarrow \mathbb{C}_{\ell}$ is a $\mathbb{Z}_{\ell}$-algebra morphism, and therefore it induces a ring morphism 
$$M_{g \times g}(\mathbb{Z}_{\ell}\llbracket T \rrbracket) \longrightarrow M_{g \times g}(\mathbb{C}_{\ell}), $$
which we denote by the same symbol.  Now, if $n \in \mathbb{N}$ and $\psi \in G(n)^{\vee}$, then we have
\begin{equation*}
\begin{aligned}
Q(t_{\tilde{\psi}}) &= {\rm ev}_{t_{\tilde{\psi}}}(Q(T)) \\
&= {\rm ev}_{t_{\tilde{\psi}}}({\rm det}(D - A_{\rho})) \\
&= {\rm det}(D - {\rm ev}_{t_{\tilde{\psi}}}(A_{\rho})).
\end{aligned}
\end{equation*}
By (\ref{imp_equal}), we have
\begin{equation*}
\begin{aligned}
{\rm ev}_{t_{\tilde{\psi}}}(A_{\rho}) &= \left(\sum_{\substack{s \in S \\ {\rm inc}(s) = (v_{i},v_{j})}}{\rm ev}_{t_{\tilde{\psi}}} \circ \rho(\alpha(s)) + \sum_{\substack{s \in S \\ {\rm inc}(s) = (v_{j},v_{i})}}{\rm ev}_{t_{\tilde{\psi}}} \circ \rho(-\alpha(s)) \right) \\
&= \left(\sum_{\substack{s \in S \\ {\rm inc}(s) = (v_{i},v_{j})}}\tilde{\psi}(\alpha(s)) + \sum_{\substack{s \in S \\ {\rm inc}(s) = (v_{j},v_{i})}}\tilde{\psi}(-\alpha(s)) \right) \\
&= \left(\sum_{\substack{s \in S \\ {\rm inc}(s) = (v_{i},v_{j})}}\psi(\alpha_{n}(s)) + \sum_{\substack{s \in S \\ {\rm inc}(s) = (v_{j},v_{i})}}\psi(-\alpha_{n}(s)) \right) \\
&= A_{\psi},
\end{aligned}
\end{equation*}
where $A_{\psi}$ is the twisted adjacency matrix, and the last equality is true by \cite[Corollary 5.2]{McGown/Vallieres:2022a}.  It follows from (\ref{3_term}) that
\begin{equation*}
\begin{aligned}
Q(t_{\tilde{\psi}}) &= {\rm det}(D - A_{\psi}) \\
&= h_{X}(1,\psi),
\end{aligned}
\end{equation*}
and this is what we wanted to show.
\end{proof}
We can now prove our main result.
\begin{theorem} \label{goal}
Let $X$ be as above and let $\alpha:S \longrightarrow \mathbb{Z}_{\ell}^{d}$ be a function for which all the graphs $X(G(n),S,\alpha_{n})$ are connected.  Consider the $\mathbb{Z}_{\ell}^{d}$-tower of graphs
$$X \longleftarrow X(G(1),S,\alpha_{1}) \longleftarrow X(G(2),S,\alpha_{2}) \longleftarrow \ldots \longleftarrow X(G(n),S,\alpha_{n}) \longleftarrow \ldots $$
and let $\kappa_{n}$ be the number of spanning trees of $X_{n} = X(G(n),S,\alpha_{n})$.  Then, there exists
$$P(X,Y) \in \mathbb{Q}[X,Y] $$
of total degree at most $d$ and of degree at most $1$ in $Y$ such that
$${\rm ord}_{\ell}(\kappa_{n}) = P(\ell^{n},n),$$
when $n$ is large enough.
\end{theorem}
\begin{proof}
Equation (\ref{imp_formula}) applied to the cover $X_{n}/X$ gives
$$\ell^{dn} \cdot \kappa_{n} = \kappa_{X} \prod_{\psi \neq \psi_{0}}h_{X}(1,\psi),$$
where the product is over all non-trivial characters of ${\rm Gal}(X_{n}/X) \simeq (\mathbb{Z}/\ell^{n}\mathbb{Z})^{d}$.  We now let
$$W^{*} = \mu_{\ell^{n}}^{d} \smallsetminus \{(1,1,\ldots,1) \}, $$
and if $\xi = (\xi_{1},\ldots,\xi_{d}) \in W^{*}$, then we set 
$$1 - \xi  = (1-\xi_{1},\ldots,1 -\xi_{d}) \in D.$$
By Theorem \ref{main_thm}, we have
\begin{equation*}
\begin{aligned}
{\rm ord}_{\ell}(\kappa_{n}) &= -dn + {\rm ord}_{\ell}(\kappa_{X}) + \sum_{\psi \neq \psi_{0}}v_{\ell}(h_{X}(1,\psi)) \\
&= -dn + {\rm ord}_{\ell}(\kappa_{X}) + \sum_{\psi \neq \psi_{0}}v_{\ell}(Q(t_{\tilde{\psi}})) \\
&= -dn + {\rm ord}_{\ell}(\kappa_{X}) + \sum_{\xi \in W^{*}}v_{\ell}(Q(1-\xi))
\end{aligned}
\end{equation*}
where $Q(T)$ is the power series in the statement of Theorem \ref{main_thm}.  Now, \cite[Theorem 5.6]{Monsky:1981} shows that there exists $E(X,Y) \in \mathbb{Q}[X,Y]$, of total degree at most $d$ and of degree at most $1$ in $Y$, such that
$$\sum_{\xi \in W^{*}}v_{\ell}(Q(1-\xi)) = E(\ell^{n},n),$$
when $n$ is large enough.  It then suffices to set
$$P(X,Y) = E(X,Y) - dY + {\rm ord}_{\ell}(\kappa_{X}) \in \mathbb{Q}[X,Y]$$
to conclude the proof.  Combining with \S \ref{towers} also proves Theorem \ref{thmA} from the introduction.
\end{proof}

Let us make a few remarks.  It is known that the coefficients of $X^{d}$ and of $Y \cdot X^{(d-1)}$ are nonnegative integers.  See \cite[Remark 2]{Monsky:1981}.  Furthermore, when $d=1$, Theorem \ref{goal} reduces to \cite[Theorem 6.1]{McGown/Vallieres:2022a} with the difference that the Iwasawa invariants $\mu,\lambda$, and $\nu$ are more easily calculated from the power series $Q(T)$ in one variable than the Greenberg coefficients when $d \ge 2$.

\section{Examples} \label{examples}
In this section, we present a few numerical examples of $\mathbb{Z}_{\ell}^{2}$-towers of graphs in the situation where the base graph $X$ is a bouquet and the function $\alpha:S \longrightarrow \mathbb{Z}_{\ell}^{2}$ takes values in $\mathbb{Z}^{2}$.  In \cite{Vallieres:2021}, \cite{McGown/Vallieres:2022}, and \cite{McGown/Vallieres:2022a}, we could calculate the Iwasawa invariants $\mu, \lambda$, and $\nu$ precisely, since we knew how far up a $\mathbb{Z}_{\ell}$-tower one had to go for the formula 
$${\rm ord}_{\ell}(\kappa_{n}) = \mu\ell^{n} + \lambda n + \nu$$ 
to hold true, and with enough $\ell$-adic precision, we could calculate $\mu$ and $\lambda$ from the power series $Q(T)$ arising in Theorem \ref{main_thm}.  Here, for each example we find some rational numbers $a,b,c,d,e \in \mathbb{Q}$ satisfying 
$${\rm ord}_{\ell}(\kappa_{n}) = a \ell^{2n} + b  n  \ell^{n} + c  \ell^{n} + d n + e$$
for a few layers only, but we have not tried to prove that those numbers are the Greenberg coefficients.  We now explain how we found candidates $a,b,c,d,e \in \mathbb{Q}$ numerically.

For this calculation, we work in $\overline{\mathbb{Q}} \subseteq \mathbb{C}$, so we do not embed everything in $\overline{\mathbb{Q}}_{\ell}$ via the embedding (\ref{embed}).  The absolute Galois group $G_{\mathbb{Q}} = {\rm Gal}(\overline{\mathbb{Q}}/\mathbb{Q})$ acts in the usual way on $G(n)^{\vee}$.  The set of orbits for this group action will be denoted by $\widehat{G(n)}(\mathbb{Q})$, and an orbit will be denoted by $\Psi$ or a similar notation.  The orbit consisting of the trivial character only will be denoted by $\Psi_{0}$.  Equation (\ref{imp_formula}) becomes at level $n$ of a $\mathbb{Z}_{\ell}^{2}$-tower over a bouquet the following formula
\begin{equation} \label{orb_prod}
\ell^{2n} \cdot \kappa_{n} = \prod_{\substack{\Psi \in \widehat{G(n)}(\mathbb{Q}) \\ \Psi \neq \Psi_{0}}}h_{X}(1,\Psi),
\end{equation}
where
$$h_{X}(1,\Psi) = \prod_{\psi \in \Psi}h_{X}(1,\psi) \in \mathbb{Z}.$$
Then, we proceed in calculating each of the integers $h_{X}(1,\Psi)$ as follows.  One has a non-canonical isomorphism $\gamma_{n}:G(n) \longrightarrow G(n)^{\vee}$ given by $(\bar{a}_{1},\bar{a}_{2}) \mapsto \psi_{(\bar{a}_{1},\bar{a}_{2})}$, where $\psi_{(\bar{a}_{1},\bar{a}_{2})}$ is defined via
$$\psi_{(\bar{a}_{1},\bar{a}_{2})}(\bar{b}_{1},\bar{b}_{2}) = \zeta_{\ell^{n}}^{\bar{a}_{1} \cdot \bar{b}_{1} + \bar{a}_{2} \cdot \bar{b}_{2}}.$$ 
The group $(\mathbb{Z}/\ell^{n}\mathbb{Z})^{\times}$ acts via the diagonal action on $G(n)$ and as such the group isomorphism $\gamma_{n}$ is equivariant.  Therefore, the orbits of the action of $(\mathbb{Z}/\ell^{n}\mathbb{Z})^{\times}$ on $G(n)$ are in one-to-one correspondence with $\widehat{G(n)}(\mathbb{Q})$.  From now on, for a positive integer $m$, we let $\varepsilon_{m}:\mathbb{Z} \longrightarrow \mathbb{R} \subseteq \mathbb{C}$ be defined via
$$a \mapsto \varepsilon_{m}(a) = (1 - \zeta_{m}^{a})(1 - \zeta_{m}^{-a}). $$
If $\alpha(s) = (b_{1,s},b_{2,s}) \in \mathbb{Z}^{2}$ for $s \in S$, \cite[Equation 17]{Vallieres:2021} implies
\begin{equation} \label{concrete_special}
h_{X}(1,\psi_{(\bar{a}_{1},\bar{a}_{2})}) = \sum_{s \in S} \varepsilon_{\ell^{n}}(\bar{a}_{1} \cdot \bar{b}_{1,s} + \bar{a}_{2} \cdot \bar{b}_{2,s}).
\end{equation}
We then calculate the orbits $\mathcal{O}$ of the action of $(\mathbb{Z}/\ell^{n}\mathbb{Z})^{\times}$ on $G(n)$, and using (\ref{concrete_special}), we calculate $h_{X}(1,\Psi)$ to high precision enough in order to be able to recognize it as an integer.  Here $\Psi$ is the orbit in $\widehat{G(n)}(\mathbb{Q})$ corresponding to $\mathcal{O}$.  Once we have all of the integers $h_{X}(1,\Psi)$, we calculate their $\ell$-adic valuation and we sum them up.  Subtracting $2n$ from this sum allows us to calculate ${\rm ord}_{\ell}(\kappa_{n})$ numerically by (\ref{orb_prod}) above.  We then go ahead and solve a system of linear equations of the form
\begin{equation*}
\begin{pmatrix}
\ell^{2(n-4)} & (n-4)\ell^{n-4} & \ell^{n-4} & n-4 & 1 & \bigm| & \text{ord}_{\ell}(\kappa_{n-4})\\
\ell^{2(n-3)} & (n-3)\ell^{n-3} & \ell^{n-3} & n-3 & 1 & \bigm| & \text{ord}_{\ell}(\kappa_{n-3})\\
\ell^{2(n-2)} & (n-2)\ell^{n-2} & \ell^{n-2} & n-2 & 1 & \bigm| & \text{ord}_{\ell}(\kappa_{n-2})\\
\ell^{2(n-1)} & (n-1)\ell^{n-1} & \ell^{n-1} & n-1 & 1 & \bigm| & \text{ord}_{\ell}(\kappa_{n-1})\\
\ell^{2n} & n\ell^n & \ell^n & n & 1 & \bigm| & \text{ord}_{\ell}(\kappa_n)
\end{pmatrix}
\end{equation*}
and this is how we get candidates $a,b,c,d,e \in \mathbb{Q}$ for each of the examples below.

\begin{enumerate}
\item Let $\ell = 2$, $S = \{s_1, s_2\}$ and $\alpha: S \longrightarrow \mathbb{Z}^2_2$ be defined via $\alpha(s_1) = (1,0)$ and $\alpha(s_2) = (0,1)$.  Then we get:
$$
\parbox{.22\linewidth}{
    \centering
        \includegraphics[scale=0.25]{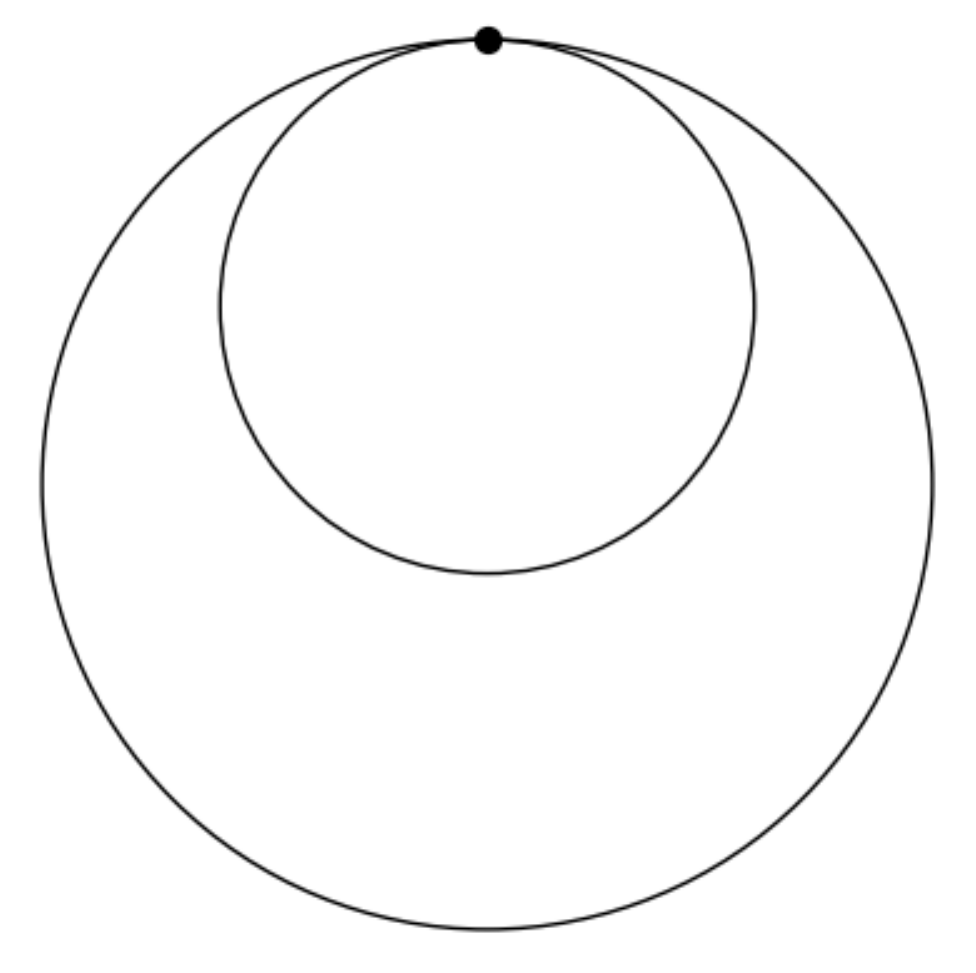}
}
\parbox{.05\linewidth}{
    $\longleftarrow$
}
\parbox{.22\linewidth}{
    \centering
        \includegraphics[scale=0.25]{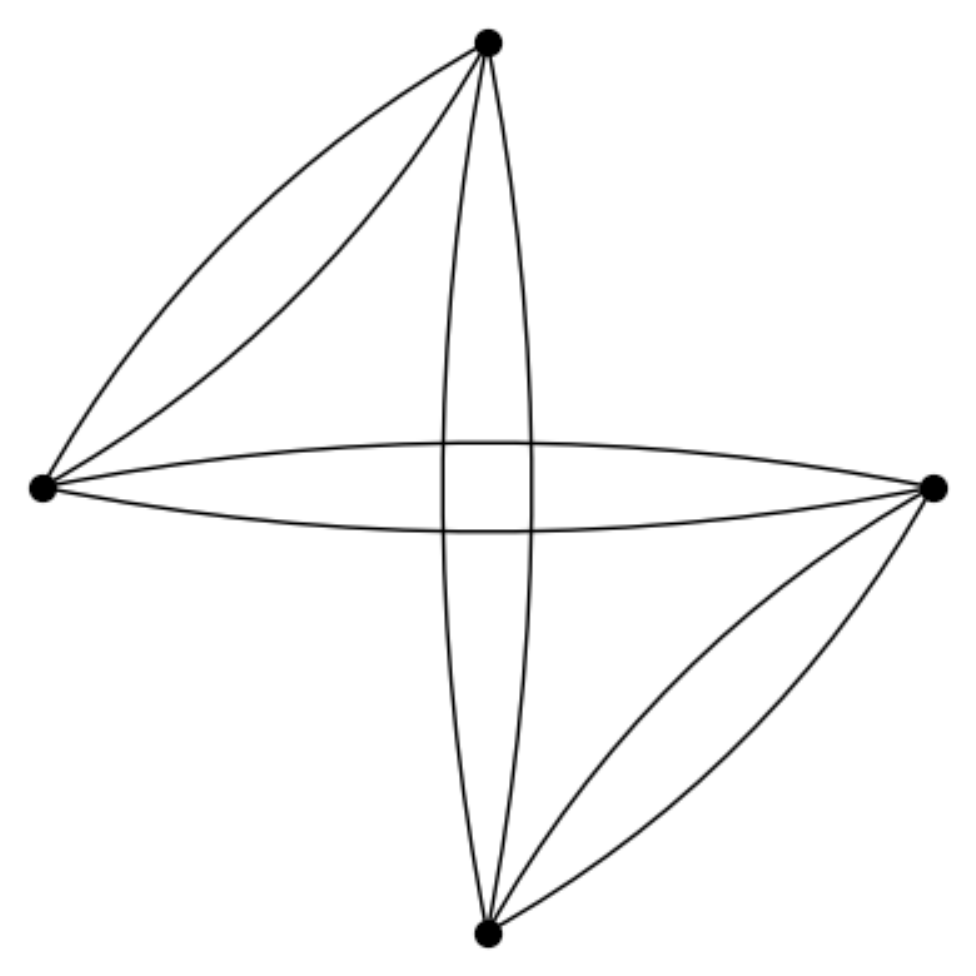}
}
\parbox{.05\linewidth}{
    $\longleftarrow$
}
\parbox{.22\linewidth}{
    \centering
        \includegraphics[scale=0.25]{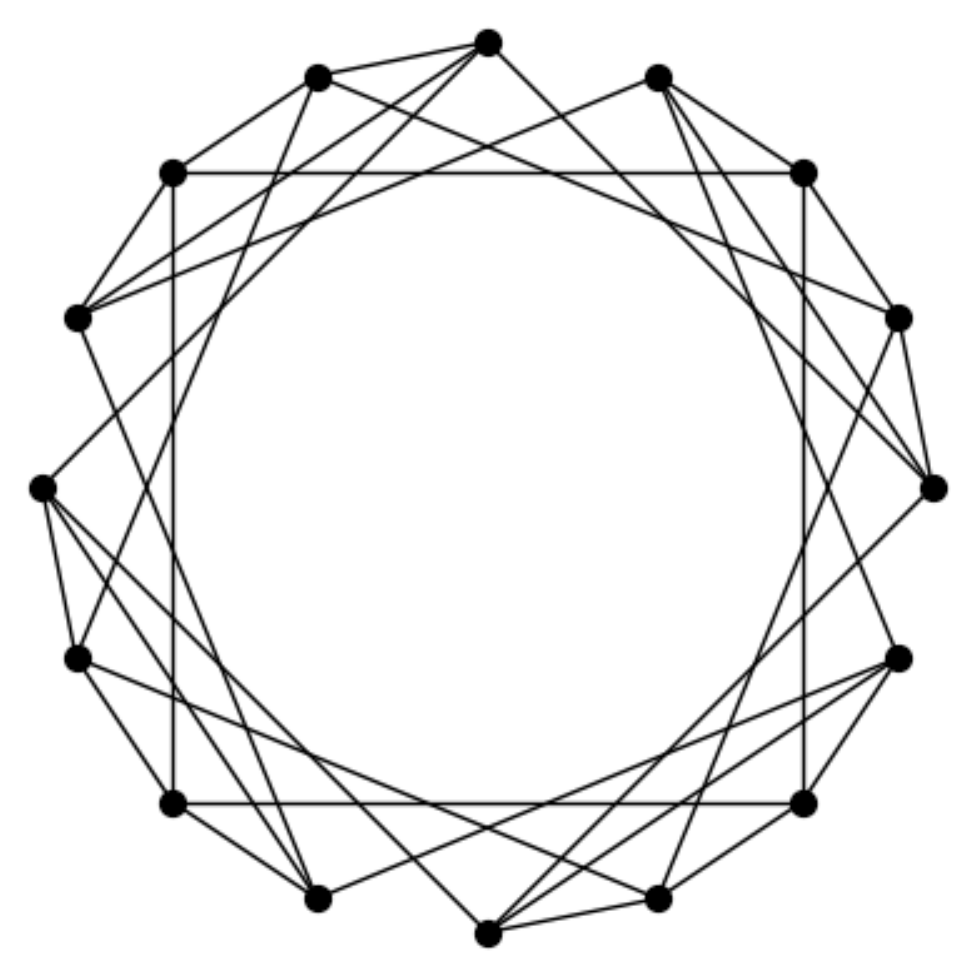}
}
\parbox{.05\linewidth}{
    $\longleftarrow$
}
\parbox{.05\linewidth}{
    $\cdots$
}
$$
We apply our algorithm in SageMath (\cite{SAGE}) to find
$$
\begin{tabularx}{0.33\textwidth} { 
    >{\centering\arraybackslash}X 
  | >{\centering\arraybackslash}X }
$n$ & $\text{ord}_2(\kappa_n)$  \\
\hline
1 & 5  \\
2 & 19  \\
3 & 61  \\
4 & 167  \\
5 & 417  \\
6 & 987  \\
7 & 2261  \\
8 & 5071  \\
9 & 11209  \\
10 & 24515 \\
\end{tabularx}
$$
which follows the pattern 
$$\text{ord}_2(\kappa_n) = 2\cdot n\cdot 2^n + 4\cdot 2^n - 6 \cdot n - 1,$$
when $1 \le n \le10$.

\item Let $\ell = 2$, $S = \{s_1, s_2, s_3, s_4\}$ and $\alpha: S \longrightarrow \mathbb{Z}^2_2$ be defined via $\alpha(s_1) = (1,0)$, $\alpha(s_2) = (1,0)$, $\alpha(s_2) = (0,1)$ and $\alpha(s_4) = (0,1)$. Then we get:
$$
\parbox{.22\linewidth}{
    \centering
        \includegraphics[scale=0.25]{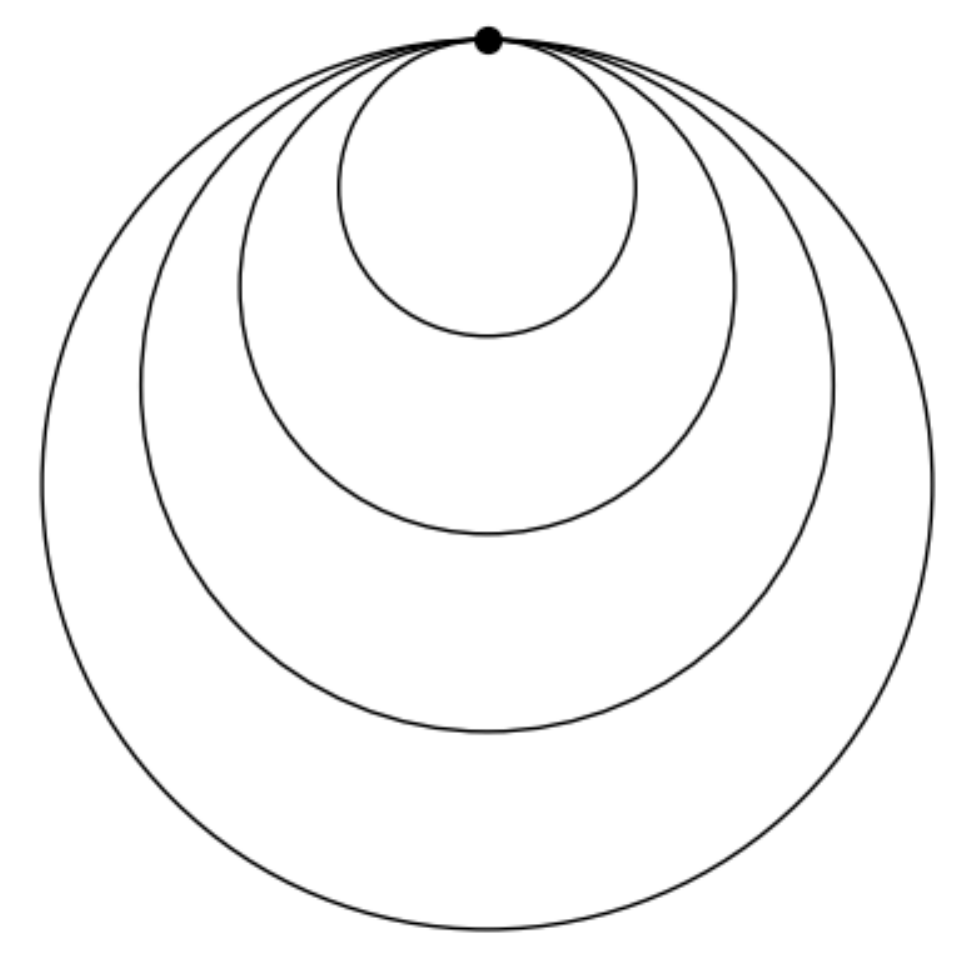}
}
\parbox{.05\linewidth}{
    $\longleftarrow$
}
\parbox{.22\linewidth}{
    \centering
        \includegraphics[scale=0.25]{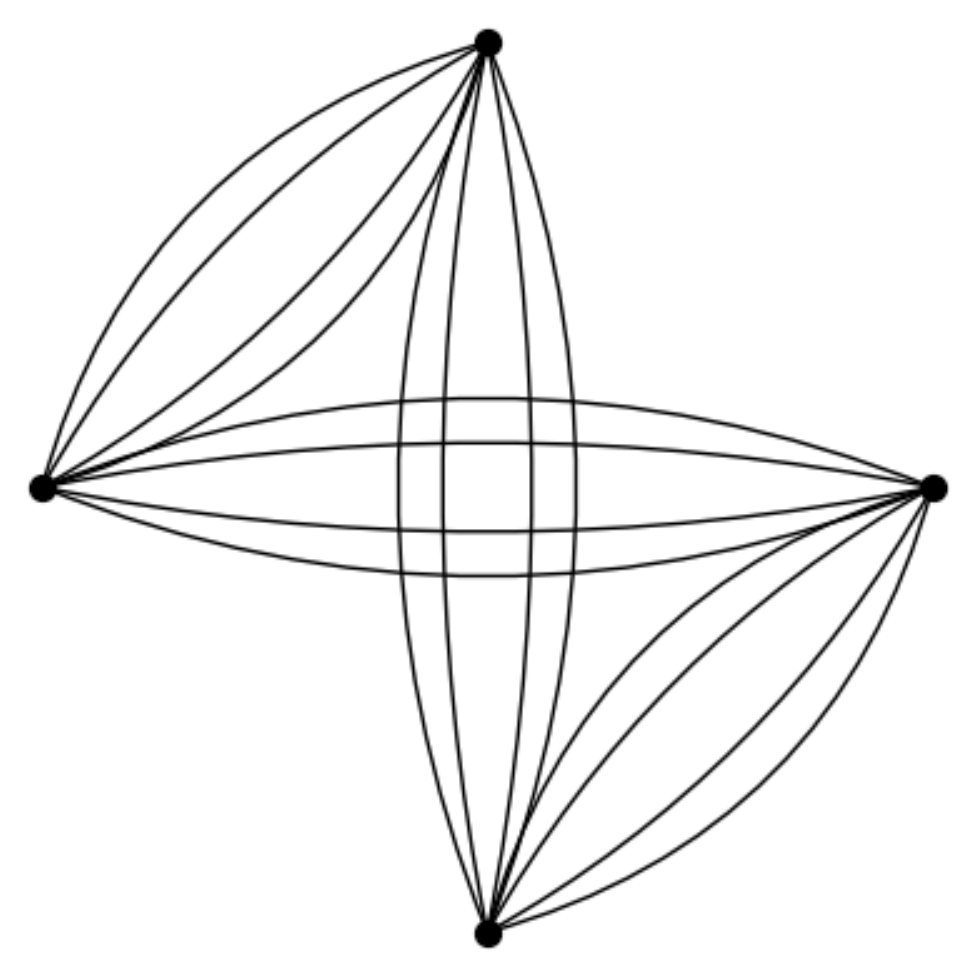}
}
\parbox{.05\linewidth}{
    $\longleftarrow$
}
\parbox{.22\linewidth}{
    \centering
        \includegraphics[scale=0.25]{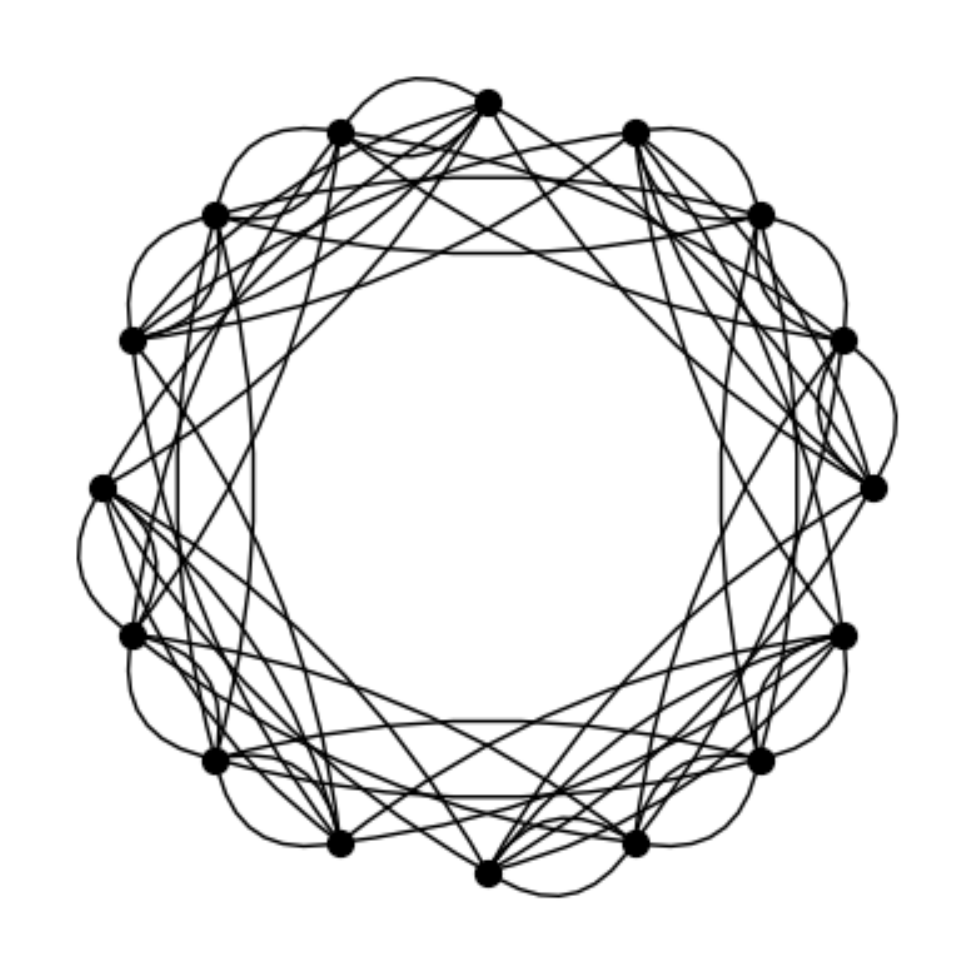}
}
\parbox{.05\linewidth}{
    $\longleftarrow$
}
\parbox{.05\linewidth}{
    $\cdots$
}
$$
Applying our algorithm reveals
$$
\begin{tabularx}{0.33\textwidth} { 
    >{\centering\arraybackslash}X 
  | >{\centering\arraybackslash}X }
$n$ & $\text{ord}_2(\kappa_n)$  \\
\hline
1 & 8  \\
2 & 34  \\
3 & 124  \\
4 & 422  \\
5 & 1440  \\
6 & 5082 \\
7 & 18644  \\
8 & 70606  \\
9 & 273352  \\
10 & 1073090\\
\end{tabularx}
$$
This time we have
$$\text{ord}_2(\kappa_n) = 2^{2n} + 2\cdot n \cdot 2^n + 4\cdot2^n - 6n - 2,$$ 
when $1 \le n \le 10$. 

\item Let $\ell = 2$, $S = \{s_1, s_2, s_3,s_4\}$ and $\alpha: S \longrightarrow \mathbb{Z}^2_2$ be defined via $\alpha(s_1) = (1,5)$, $\alpha(s_2) = (0,3)$, $\alpha(s_3) = (1,2)$, and $\alpha(s_4) = (0,1)$.  Then we get:
$$
\parbox{.22\linewidth}{
    \centering
        \includegraphics[scale=0.25]{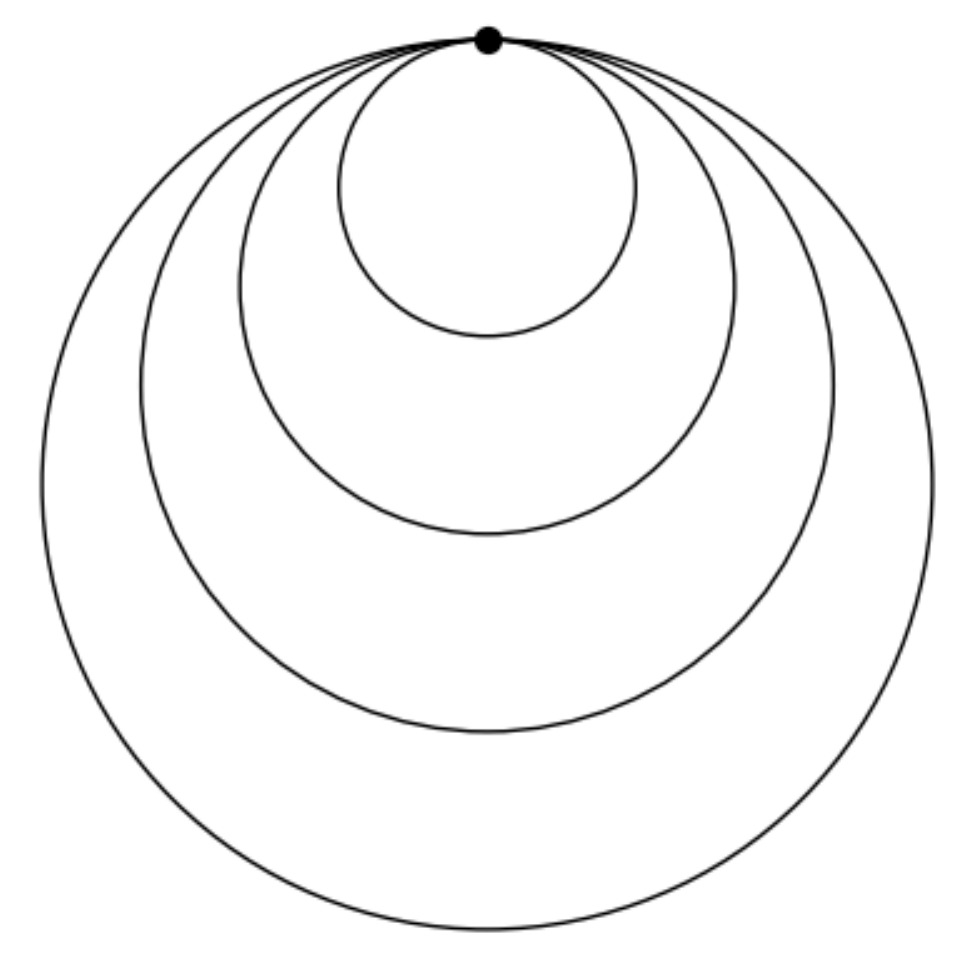}
}
\parbox{.05\linewidth}{
    $\longleftarrow$
}
\parbox{.22\linewidth}{
    \centering
        \includegraphics[scale=0.25]{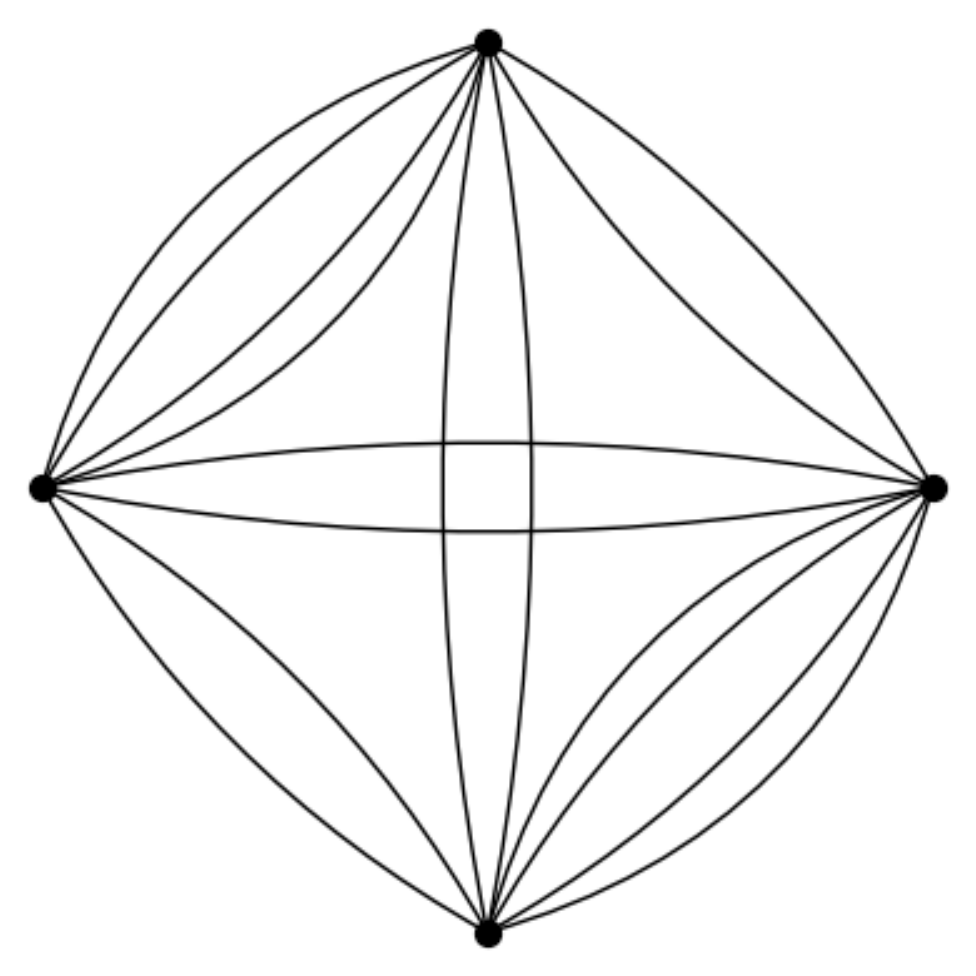}
}
\parbox{.05\linewidth}{
    $\longleftarrow$
}
\parbox{.22\linewidth}{
    \centering
        \includegraphics[scale=0.25]{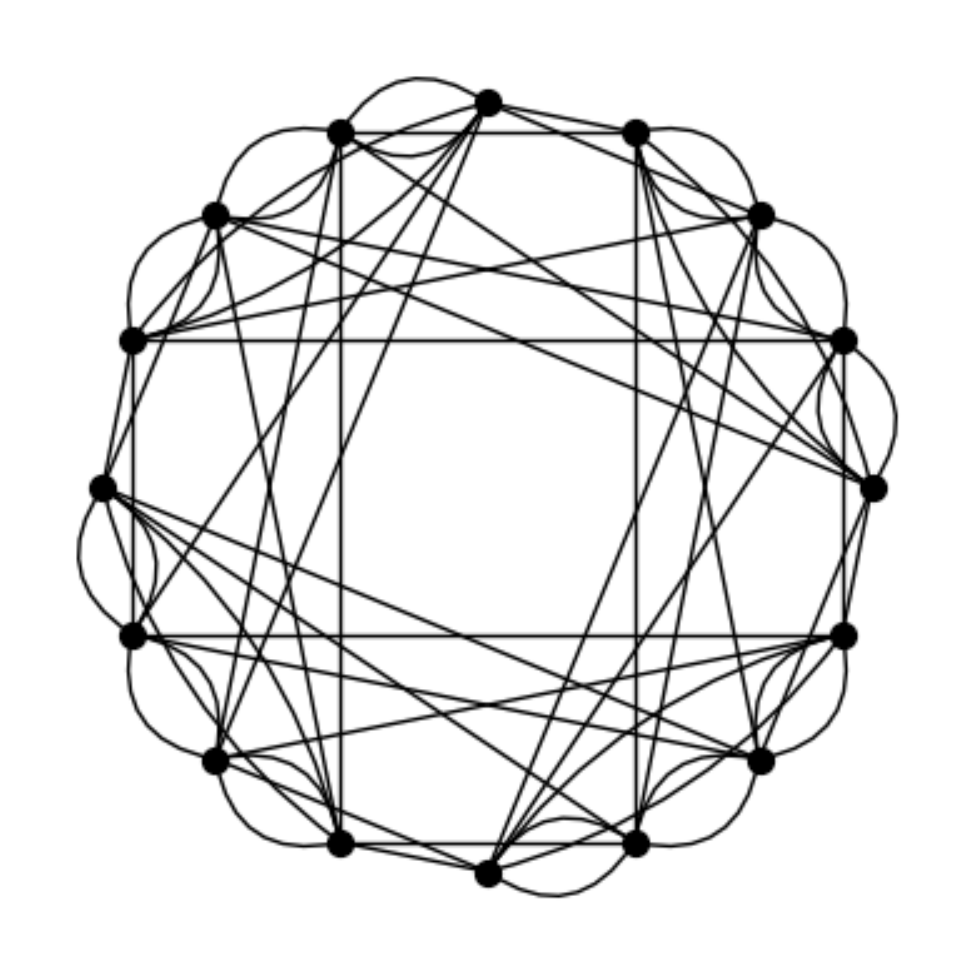}
}
\parbox{.05\linewidth}{
    $\longleftarrow$
}
\parbox{.05\linewidth}{
    $\cdots$
}
$$
After an application of our algorithm we get the data
$$
\begin{tabularx}{0.33\textwidth} { 
    >{\centering\arraybackslash}X 
  | >{\centering\arraybackslash}X }
$n$ & $\text{ord}_2(\kappa_n)$  \\
\hline
1 & 5  \\
2 & 19  \\
3 & 65  \\
4 & 179  \\
5 & 403 \\
6 & 887 \\
7 & 1923  \\
8 & 4127 \\ 
9 & 8795 \\
10 & 18647 \\
\end{tabularx}
$$
from which we obtain the formula 
$$\text{ord}_{2}(\kappa_n) = n\cdot2^n + \frac{33}{4} \cdot 2^n - 4n - 1,$$ 
for $4 \le n \le 10$.

\item Let $\ell = 3$, $S = \{s_1, s_2\}$ and $\alpha: S \longrightarrow \mathbb{Z}^2_3$ be defined via $\alpha(s_1) = (1,0)$ and $\alpha(s_1) = (0,1)$. Then we get:
$$
\parbox{.22\linewidth}{
    \centering
        \includegraphics[scale=0.25]{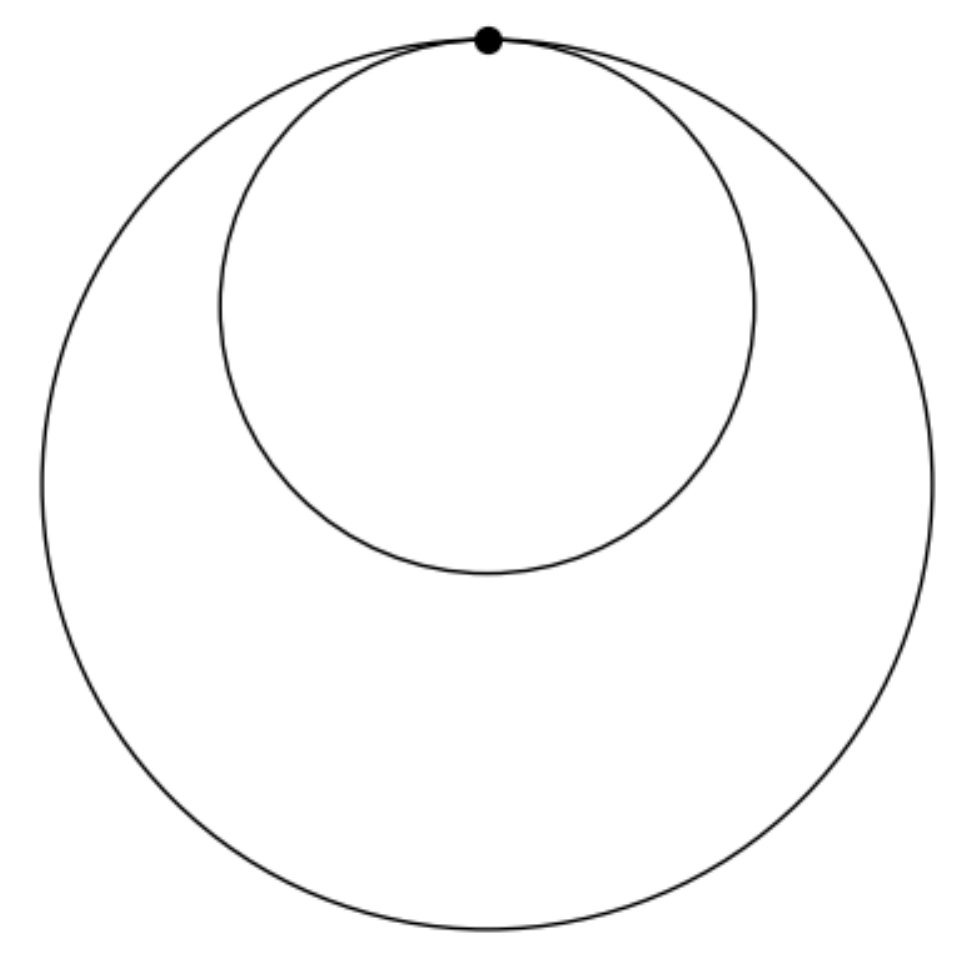}
}
\parbox{.05\linewidth}{
    $\longleftarrow$
}
\parbox{.22\linewidth}{
    \centering
        \includegraphics[scale=0.25]{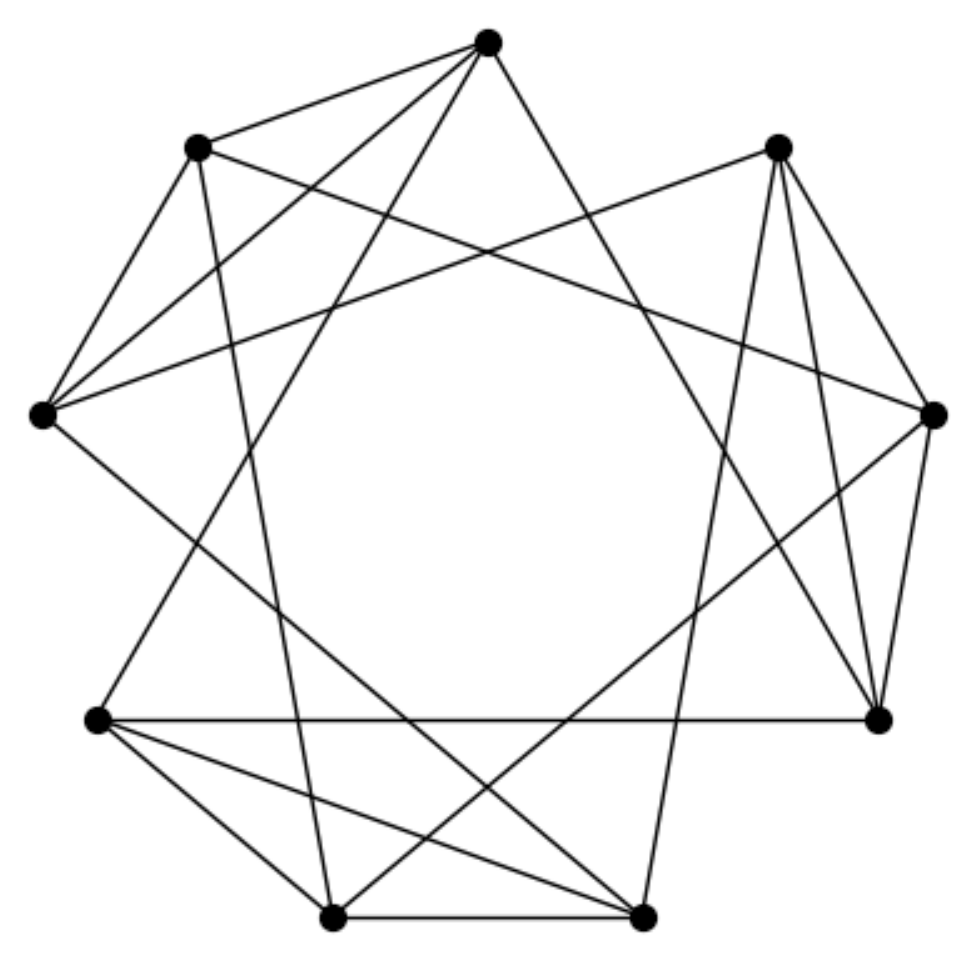}
}
\parbox{.05\linewidth}{
    $\longleftarrow$
}
\parbox{.22\linewidth}{
    \centering
        \includegraphics[scale=0.25]{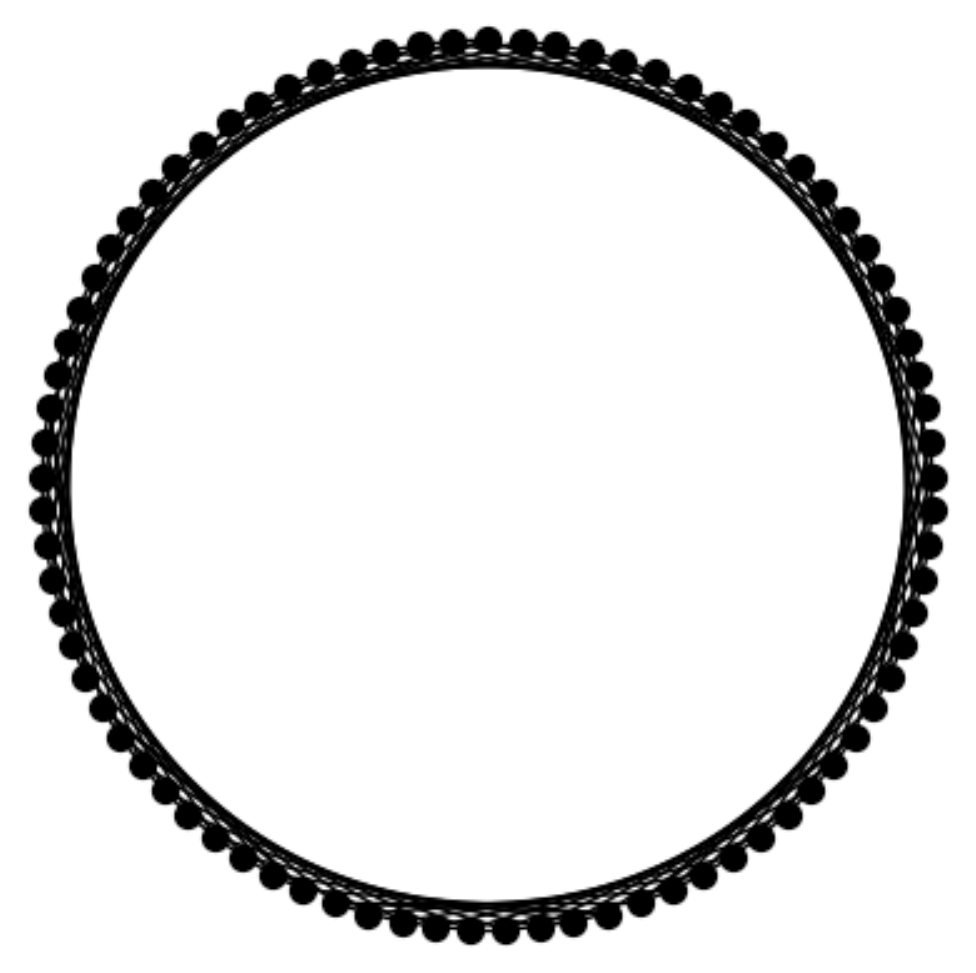}
}
\parbox{.05\linewidth}{
    $\longleftarrow$
}
\parbox{.05\linewidth}{
    $\cdots$
}
$$
After applying our algorithm we get the table of values
$$
\begin{tabularx}{0.33\textwidth} { 
    >{\centering\arraybackslash}X 
  | >{\centering\arraybackslash}X }
$n$ & $\text{ord}_3(\kappa_n)$  \\
\hline
1 & 6  \\
2 & 28  \\
3 & 98  \\
4 & 312  \\
5 & 958  \\
6 & 2900 \\
7 & 8730 \\
\end{tabularx}
$$
from which we obtain the following formula 
$$\text{ord}_3(\kappa_n) = 4\cdot 3^n - 2n - 4,$$ 
for $1 \le n \le 7$.

\item Let $\ell = 3$, $S = \{s_1, s_2, s_3\}$ and $\alpha: S \longrightarrow \mathbb{Z}^2_3$ be defined via $\alpha(s_1) = (1,0)$, $\alpha(s_2) = (2,3)$, and $\alpha(s_4) = (1,1)$. Then we get:
$$
\parbox{.22\linewidth}{
    \centering
        \includegraphics[scale=0.25]{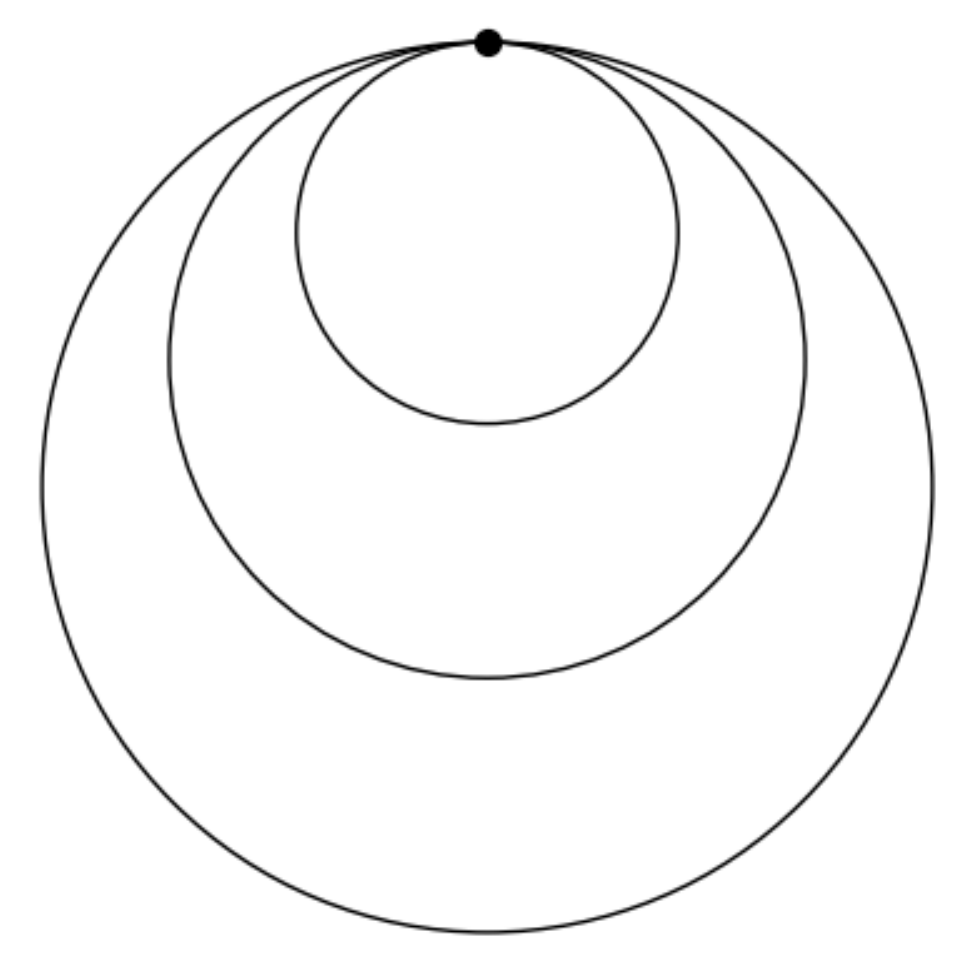}
}
\parbox{.05\linewidth}{
    $\longleftarrow$
}
\parbox{.22\linewidth}{
    \centering
        \includegraphics[scale=0.25]{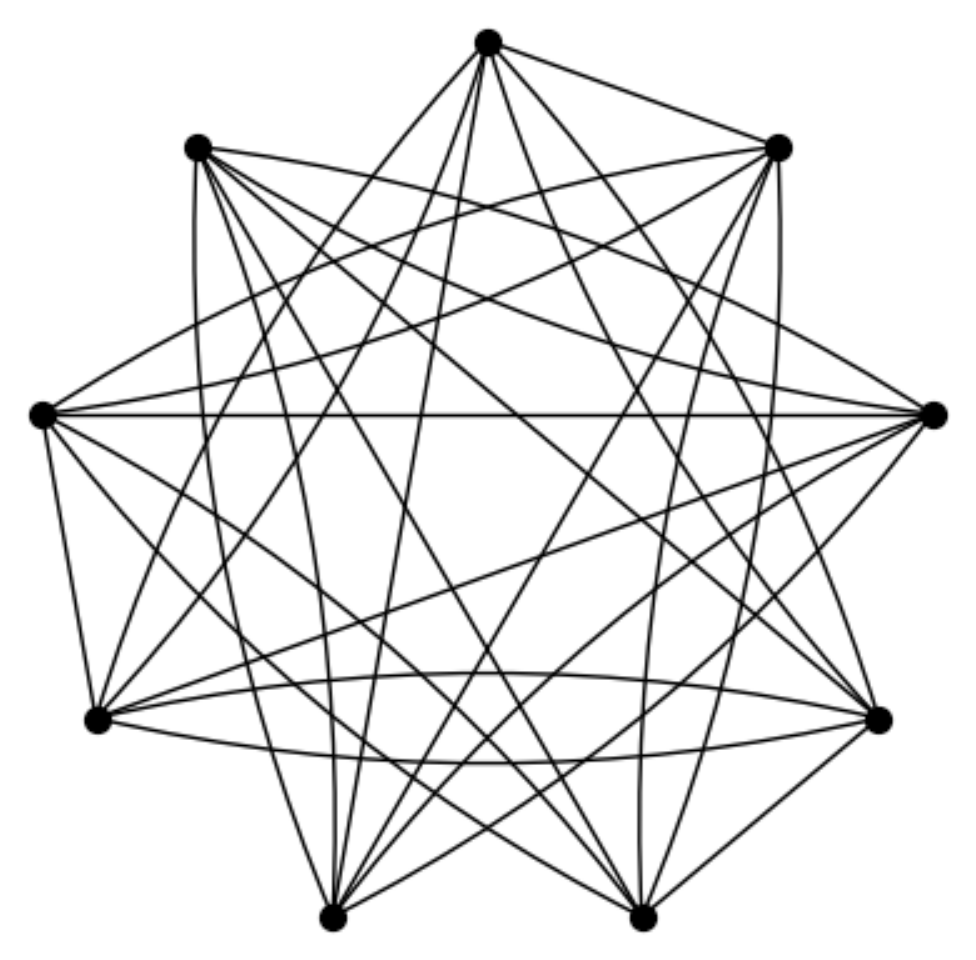}
}
\parbox{.05\linewidth}{
    $\longleftarrow$
}
\parbox{.22\linewidth}{
    \centering
        \includegraphics[scale=0.25]{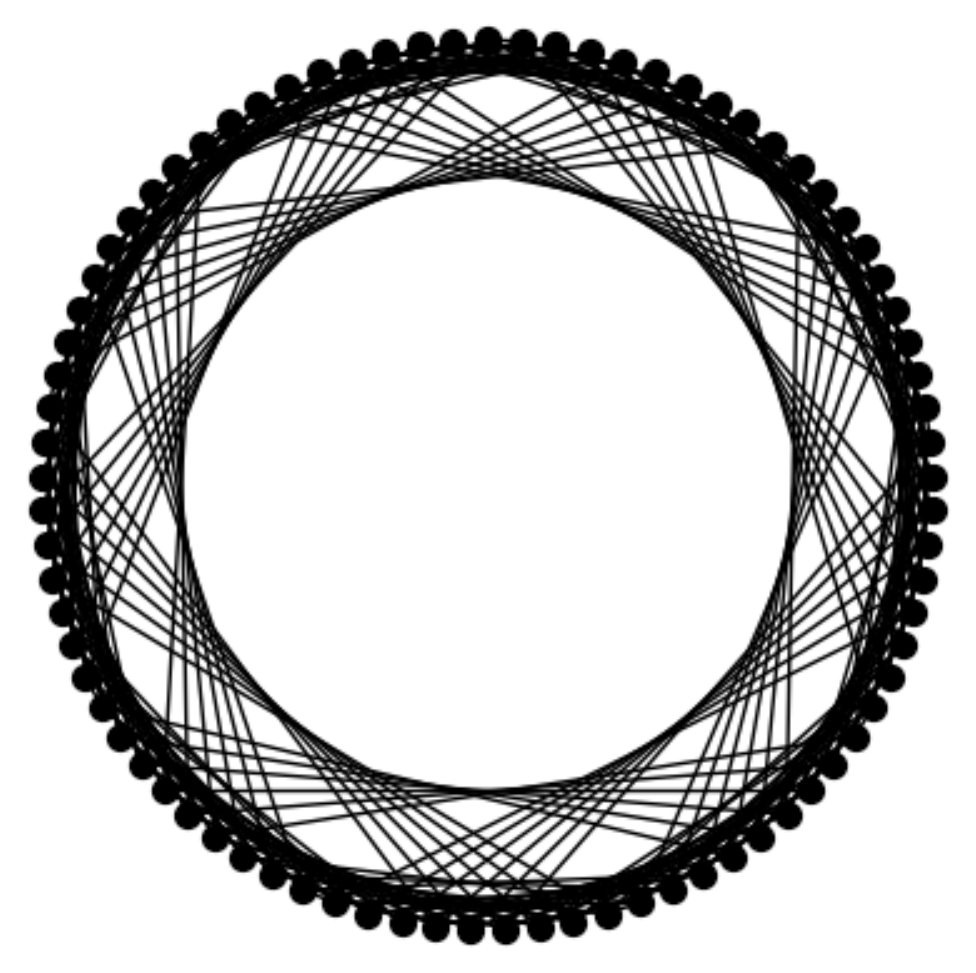}
}
\parbox{.05\linewidth}{
    $\longleftarrow$
}
\parbox{.05\linewidth}{
    $\cdots$
}
$$
and after an application of our algorithm we get
$$
\begin{tabularx}{0.33\textwidth} { 
    >{\centering\arraybackslash}X 
  | >{\centering\arraybackslash}X }
$n$ & $\text{ord}_3(\kappa_n)$  \\
\hline
1 & 10  \\
2 & 48  \\
3 & 166  \\
4 & 524  \\
5 & 1602  \\
6 & 4840 \\
7 & 14558  \\
\end{tabularx}
$$
from which we obtain the formula 
$$\text{ord}_3(\kappa_n) = \frac{20}{3}3^n - 2n - 8,$$ 
for $1 \le n \le 7$.

\end{enumerate}

\bibliographystyle{plain}
\bibliography{main}

\end{document}